%% file: main.tex
\documentclass[10pt]{amsart}
\usepackage[margin=2cm]{geometry}
\usepackage{import}
\usepackage{graphics, amsmath,amssymb, enumitem, comment, url, mathtools}
\usepackage{graphicx}
\usepackage{epsfig}
\usepackage{tikz}
\usepackage{forest}
\usepackage{bm}
\usetikzlibrary{trees}
\usepackage{enumitem}

\usepackage[titles]{tocloft}
\newlength\mylength

\usepackage{graphicx}

\usepackage{subcaption}
\usepackage{titlesec}

\titleformat{\section}
{\large\filcenter\scshape}
{\thesection.}{0.5em}{}

\titleformat{\subsection}[runin]
{\normalfont\bfseries}
{\thesubsection.}{.5em}{}[. ]
\titlespacing{\subsection}
{\parindent}{1ex}{5pt}

\newif\ifcolorcomments

\newcommand{\blalpha}{\boldsymbol{\alpha}}

\newcommand{\knm}{\boldsymbol{k}}
\newcommand{\blv}{\boldsymbol{v}}
\newcommand{\blu}{\boldsymbol{u}}
\newcommand{\blw}{\boldsymbol{w}}
\newcommand{\bbeta}{\boldsymbol{\beta}}

\newcommand{\func}[1]{\psi_{\alpha_{#1}}}
\newcommand{\Z}{\mathbb{Z}}

\usepackage{graphicx,xcolor}
\graphicspath{{figures/}} % путь, где искать картинки
\newcommand{\executeiffilenewer}[3]{%
	\ifnum\pdfstrcmp{\pdffilemoddate{#1}}%
	{\pdffilemoddate{#2}} > 0 {\immediate\write18{#3}}\fi}

\colorcommentstrue
\usepackage{amssymb}
\usepackage{amsmath,amsthm}
\usepackage{hyperref}

\usepackage{colortbl}

\newtheorem{theorem}{Theorem}[section]
\newtheorem{lemma}[theorem]{Lemma}
\newtheorem{proposition}[theorem]{Proposition}

\theoremstyle{definition}

\newtheorem{remark}[theorem]{Remark}
\newtheorem{example}[theorem]{Example}

\usepackage{graphicx}

\makeatletter
\DeclareRobustCommand{\cev}[1]{%
	{\mathpalette\do@cev{#1}}%
}
\newcommand{\do@cev}[2]{%
	\vbox{\offinterlineskip
		\sbox\z@{$\m@th#1 x$}%
		\ialign{##\cr
			\hidewidth\reflectbox{$\m@th#1\vec{}\mkern4mu$}\hidewidth\cr
			\noalign{\kern-\ht\z@}
			$\m@th#1#2$\cr
		}%
	}%
}
\makeatother

\newcommand{\N}{\mathbb N}% by default, $\N$ is American naturals $\{1,2,\ldots\}$
\newcommand{\Q}{\mathbb Q}

\newcommand{\R}{\mathbb R}

\renewcommand{\Z}{\mathbb Z}

\DeclareMathOperator{\Sym}{Sym}

\renewcommand{\text}{\textup}

\newcommand{\NPC}[1]{\ignorespaces}

\newif\ifdraft\drafttrue

\def\N{\mathbb N}
\def\Z{\mathbb Z}
\def\Q{\mathbb Q}

\def\R{\mathbb R}

\def\a{\boldsymbol{\alpha}}

\usepackage{calc}

\IfFileExists{marvosym.sty}{
	\RequirePackage{marvosym} 
}{
	
}

\IfFileExists{wasysym.sty}{
	\RequirePackage{wasysym}\renewcommand{\emptyset}{{\diameter}}
}{

}

\newcommand*{\myDots}{\ifmmode\mathellipsis\else.\kern-0.07em.\kern-0.07em.\fi}

\DeclareMathOperator{\pr}{Pr}

\usepackage{enumitem}

\newcommand {\ignore}[1] {}

\usepackage{mathtools}
\usepackage{ longtable}
\begin{document}
	
	\title{Permutation of values of irrationality measure functions}

	\author{Victoria Rudykh}
	%\address{Victoria Rudykh,  Faculty of Mathematics, Technion, Haifa, Israel.}
	\date{}
	
	\import{./}{abstract.tex}
	\maketitle
	\import{./}{several_functions.tex}

	\import{./}{bibliography.tex}

\end{document}

%% file: abstract.tex
\begin{abstract}
	For an irrational number $\alpha\in\mathbb{R}$ we consider its irrationality measure function $$ \psi_\alpha(t) = \min_{1\le q\le t,\, q\in\mathbb{Z}} \| q\alpha \|. $$
	Let  $\bm{\alpha} = (\alpha_1, \dots, \alpha_n)$ be $n$-tuple of pairwise independent irrational numbers.  For each $t \in \R_{>1}$ irrationality measure functions $\psi_{\alpha_1}, \dots, \psi_{\alpha_n}$ can be written in an increasing order 
	$$\psi_{\alpha_{v_1}}(t) > \psi_{\alpha_{v_2}}(t) > \dots > \psi_{\alpha_{v_{n-1}}}(t) > \psi_{\alpha_{v_n}}(t).$$
	We consider the vector of functions $\blv_{\a}(t): \R_{>1} \rightarrow S_n$ associated to this order and defined as 
	$$\blv_{\a}(t) = ( v_1, v_2, \dots, v_{n-1}, v_n ).$$
	Let $\knm(\bm{\alpha})$ be the number of infinitely occurring different values of $\blv_{\a}(t)$. It is known that if $\knm(\a)= k$ we have
	$ n \leq \frac{k(k+1)}{2}.$ 
	At the same time,  for $k \geq 3$ and $n = \frac{k(k+1)}{2}$ 
	there exists an $n$-tuple $\a$ with $\knm(\a) = k$.

	In this work we define a $k$-cyclic permutation $\pi$ and prove that in the extremal case $n = \frac{k(k+1)}{2}, \ \knm(\a) = k$ the set of successive values of $\blv_{\a}(t)$ is an orbit of  $\pi$.
\end{abstract}

%% file: several_functions.tex
\section{Introduction}
For an irrational number $\alpha\in\mathbb{R}$ we define its irrationality measure function
$$ \psi_\alpha(t) = \min_{1\le q\le t,\, q\in\mathbb{Z}} \| q\alpha \| ,$$
where $||\cdot || $ denotes the distance to the nearest integer.
It is a well-known fact that it is a piecewise constant decreasing function and
\begin{equation*}
	\psi_\alpha(t) = \| q_m\alpha \| \text{  for  } q_m\le t < q_{m+1}, 
\end{equation*}
where $q_m$ are the denominators of the convergents $\frac{p_m}{q_m} = [a_0; a_1, \dots, a_m]$ for continued fraction $\alpha = [a_0; a_1, a_2, \dots].$ 

In 2010 Kan and Moshchevitin \cite{KM} proved a surprising result about two irrationality measure functions.
\begin{theorem}[Kan, Moshchevitin, 2010]
	\label{th:KM}
	For any two  different irrational numbers $\alpha, \beta$ such that $\alpha \pm \beta \not\in \Z$ the difference function
	\begin{equation*}
		\psi_\alpha(t)  - \psi_\beta(t)
	\end{equation*}
	changes its sign infinitely many times as $t \rightarrow +\infty$
\end{theorem}
In the last decade several papers were devoted to a more detailed analysis of the mutual behavior of two irrationality measure functions (\cite{Dub17}, \cite{Mos19}, \cite{Shu22}, \cite{RS24}). The combinatorial behavior of several irrationality measure functions was recently considered in \cite{MM21}, \cite{Rud22}.

To formulate our result, we first introduce the necessary notions and notation. We are interested in the asymptotic behavior of the functions, so throughout the paper $T$ will denote a large enough positive real number. We call a real vector $\blalpha = (\alpha_1, \dots, \alpha_n) \in \R^n $ to be an $n$-tuple of pairwise independent numbers if there exists $T \in \R_{0}$ such that for $i \neq j$
\begin{equation}
	\label{def:T}
	\forall \, t \geq T \ \ \ \func{i}(t) \neq \func{j}(t).
\end{equation}
For simplification of notation we will write $\psi_i = \func{i}$.  For $t \geq T$ we rewrite the set of values 
\begin{equation*}
	\psi_{1}(t), \ \psi_{2}(t), \ \ldots \ , \ \psi_{n-1}(t), \ \psi_{n}(t)
\end{equation*}
in the increasing order
\begin{equation}
	\label{int:eq1}
	\psi_{v_1}(t) > \psi_{{v_2}}(t) > \ldots >\psi_{{v_{n-1}}}(t)> \psi_{{v_n}}(t).
\end{equation}
Now (\ref{int:eq1}) defines the \textit{vector of order of functions} 
\begin{equation*}
	\blv_{\a}: \R_{\geq T} \rightarrow S_n, \quad \blv_{\a}(t) = ( v_1, v_2, \dots, v_{n-1}, v_n ),
\end{equation*} 
where $S_n$ is the set of all permutations of elements $\{1, \dots, n\}$.
We will refer to $\blv_{\a}(t)$ only as  \textit{vector}.

By $\knm(\blalpha)$ we denote the number of different values of vector $\blv_{\a}(t)$ which occur infinitely many times when $t \rightarrow \infty$. Formally, 
\begin{equation*}
	\knm(\blalpha) = \max\{k : \exists \text{ different } \blv_1, \dots, \blv_k \in S_n: \forall \, 1 \leq j \leq k \ \ \forall \, t > 0 \ \ \exists \, t' > t \ \ \blv_{\a}(t')=\blv_j\}.
\end{equation*}
Standard metric argument shows that for almost all $\a$ (in the sense of Lebesgue measure in $\R^n$), we have $\knm(\a) = n!$.
In 2021 Manturov and Moshchevitin \cite{MM21} proved the following result about $\knm(\blalpha)$. 
\begin{theorem}
	\label{th:MM}
	For $k \geq 3$ and $n = \frac{k(k+1)}{2}$ there exists $n$-tuple $\blalpha$ of pairwise independent numbers with
	\begin{equation*}
		\knm(\alpha) = k.
	\end{equation*}
\end{theorem}
\noindent This result turned out to be optimal and in \cite{Rud22}
the author proved the following statement.
\begin{theorem}
	\label{th:upper_bound}
	The size of $n$-tuple $\blalpha = (\alpha_1, \dots, \alpha_n)$ of pairwise independent numbers with $\knm(\blalpha) = k$ necessarily satisfies the inequality
	\begin{equation*}
		n \leq \frac{k(k+1)}{2}.
	\end{equation*}
\end{theorem}
\noindent In this paper we completely describe the structure of sequences of vectors under the extremal condition ${\knm(\blalpha) = k},$ ${n = \frac{k(k+1)}{2}}$.
We will show that in this case the set of successive values $\blv_\alpha(t)$ has a very specific structure, namely, it is an orbit of a $k$-cyclic permutation $\pi$, which we define in Section \ref{sec:main result}.

\section{Main result}\label{sec:main result}
\subsection{Special enumeration and permutation $\pi$}\label{subsec:pi}
Let $n = \frac{k(k+1)}{2}$. We first introduce a new enumeration for elements $u_i$ of the vector $\blu = (u_1, \dots, u_n) \in S_n$ and then define a permutation $\pi$ on $\blu$.  The reason for this special enumeration is explained in Section \ref{sec:4}. 

A natural way of enumeration is by natural numbers $1, 2, \dots, n-1, n$. However, we enumerate them by pairs $(j, \, l), \ 1 \leq j \leq k, \ 1 \leq l \leq j$ in the following way
\begin{equation}
	\label{int:eq10}
	\blu = (\underbrace{u_{1,1}, \, u_{1,k},\, u_{1,k-1},\, \dots,\, u_{1,2}}_k, \, \underbrace{u_{2,2},\, u_{2,k}, \dots, \, u_{2,3}}_{k-1}, u_{3,3}, \dots, \, \underbrace{u_{k-1,k-1},\, u_{k-1,k},}_2\, \underbrace{u_{k,k}}_1  ).
\end{equation}
\import{./}{pic_10.tex} 
Here, the $k$ elements of the first block have indices $(1,l), \ k \geq l \geq 1$, the next $k-1$ elements of the second block have indices $(2, l), \  k-1 \geq l \geq 1$, and so on. The last block consists of just one element with the index $(k,k)$. \mbox{In the $l$th  block of length $k-l + 1$ the index of first element is $(l,l)$,} and the second indices of the following $l-1$ elements decrease from $k$ to $l+1$. More formally, 
\begin{equation*}
	%\label{int:eq2}
	u_i = \begin{cases}
		u_{j,j} & \text{ if } \ i = 1 + \sum\limits_{m = 0}^{j-2}(k - m), \  1 \leq j \leq k; \\ 
		u_{j,l}, & \text{ if } \ i = k - l + 2 + \sum\limits_{m = 0}^{j-2} (k - m), \ \ 1 \leq j \leq k-1, \ j+1 \leq l \leq k.
	\end{cases}
\end{equation*}
We represented this enumeration as a triangular $k \times k$ diagram (see the left triangle in Figure \ref{fig:pi}). In the first column we write the first $k$ elements, $u_{1,1}, u_{1,k}, \dots, u_{1,2}$, from top to bottom. In the second column we write the next $k-1$ elements $u_{2,2}, u_{2,k}, \dots, u_{2,3}$, and so on, until the last element $u_{k,k}$, which we write in the last $k$th column.

We define the permutation $\pi: (u_1, \dots, u_n) \mapsto (\pi(u_1), \dots, \pi(u_n))$ by means of Figure \ref{fig:pi} and in terms of enumeration (\ref{int:eq10}). The permutation $\pi$ acts as follows. Elements of the first row are cyclically shifted to the left, where the element $u_{1,1}$ moves to the last position in the row. The column of the elements $u_{1,k}, \dots, u_{1,2}$ is placed horizontally under the first row in reverse order. All other elements are shifted down by one row. After the shifts, the elements $u_{i,j}$ are written again as a vector according to the shifted diagram. As before, we first write elements from the first column, from top to bottom, then from the second column in the same manner, and so on, until the last element in the last column.
Formally, 
\begin{align*}
	&\pi(\blu)_{i,i} = u_{i+1,i+1}, \ \ 1 \leq i \leq k-1, \\
	\nonumber
	&\pi( \blu )_{k,k} = u_{1,1}, \\
	\nonumber
	&\pi(\blu)_{i,k} = u_{1,i + 1}, \ 1 \leq i \leq k-1 \\
	\nonumber
	&\pi (\blu)_{i,j} = u_{i+1,j+1}, \ \ 1 \leq i \leq k-2, \ i+1 \leq j \leq k-1.
\end{align*}
In Section \ref{ch_2:auxiliary lemmas} we analyze this permutation and prove that it is $k$-cyclic (Lemma \ref{lm:order}).
\import{./}{pic_7.tex} 
\begin{comment}
	\begin{equation*}
		\pi: \begin{pmatrix}
			u_{1,0} & u_{1,k} & u_{1,k-1} & \dots & u_{1,2} & u_{2,0} & u_{2,k} & \dots & u_{2,3} & u_{3,0} & u_{3,k} & \dots & u_{k-1,0} & u_{k-1,k} & u_{k,0}  \\
			u_{2,0} & u_{1,2} & u_{2,k} & \dots & u_{2,3} & u_{3,0} & u_{1,3} & u_{3,k} & \dots & u_{k-1,0} & u_{1,k-1} & u_{k-1,k} & u_{k,0} & u_{1,k} & u_{1,0}
		\end{pmatrix}
	\end{equation*}
\end{comment}
\begin{example}
	In order to better describe the action of $\pi$ we consider an example when $k = 5$ and $n = 10$ (see Figure \ref{fig:perm}). As defined before, an element $\blu = (u_1, \dots, u_{10}) \in S_{10}$ is enumerated as
	\begin{equation*}
		\blu = (u_{1,1}, \, u_{1,5},\, u_{1,4}, \,  u_{1,3}, \, u_{1,2},\, u_{2,2},\,  u_{2,5}, \,  u_{2,4}, \, u_{2,3}, \, u_{3,3}, \, u_{3,5}, \,  u_{3,4}, \,  u_{4,4}, \,  u_{4,5}, \, u_{5,5}).
	\end{equation*}
	And $\pi$ acts on it as
	\begin{align*}
		\pi(\blu) = (
		u_{2,2}, \, u_{1,2}, \, u_{2,5}, \,  u_{2,4}, \, u_{2,3}, \,  u_{3,3},\,  u_{1,3}, \, u_{3,5}, \,  u_{3,4}, \, u_{4,4}, \, u_{1,4}, \, u_{4,5}, \, u_{5,5}, \,  u_{1,5},\, u_{1,1}).
	\end{align*}
\end{example}
\subsection{Main Theorem}
Let $\blalpha$ be an $n$-tuple of $n = \frac{k(k+1)}{2}$ pairwise independent numbers with $\knm(\blalpha) =k$, $t_0 \in \R_{>1}.$ We  define inductively the sequence $\{t_i(\a, t_0)\}_{i=1}^\infty$ of moments of permutation changes by
\begin{equation}
	\label{def:ti}
	t_{i+1}(\a, t_0) = \min\{t > t_i(\a): \ \blv_{\a}(t) \neq \blv_{\a}(t_i)\}, \ i \in \N_0.
\end{equation}
We will denote this sequence only as $\{t_i\}_{i=0}^\infty$, keeping in mind that it depends on $\a$ and $t_0$.

We now formulate our main result.
\begin{theorem}
	\label{th:uniqness}
	Let $\blalpha$ be an $n$-tuple of pairwise independent irrational numbers such that $n = \frac{k(k+1)}{2}$ and $\knm(\blalpha) = k$. Then there exists $T \in \R_{>1}$ such that for all $t_0 \geq T$ and for all $i \in \N_0$ we have
	\begin{equation*}
		\blv_{\a}(t_{i+1}) = \pi(\blv_{\a}(t_i)),
	\end{equation*}
	where $\{t_{i}\}_{i=1}^{\infty} = \{t_{i}(\a, t_0)\}_{i=1}^{\infty}$ is the sequence of moments of permutation changes defined in (\ref{def:ti}).
\end{theorem}
We will prove a stronger result in Theorem \ref{th:main}, which completely describes the behavior of irrationality measure functions. Theorem \ref{th:uniqness} will follow as a consequence.

\subsection{Open questions}
This paper characterizes the critical case $k = \knm(\a), \, n = \frac{k(k+1)}{2}$. It may be interesting to see what happens when we loosen the condition on the number of functions. For example, can we describe completely the behavior of $n = 5$ functions that have $k = 3$ permutations.

Another question is for what sets of permutations can we construct an $n$-tuple of numbers $\a$ such that the vector $\blv_{\a}(t)$ takes values only from the given set.

\section{Preliminaries} \label{ch_2:auxiliary lemmas}
In this section we recall known results and prove auxiliary lemmas that we need for the proof of Theorem  \ref{th:uniqness}.
\subsection{Continued fractions} 
Throughout the rest of the paper we assume that
$\alpha\notin\mathbb{Q}$ and denote its continued fraction expansion as
\begin{equation*}
	\alpha =  a_0 + \cfrac{1}{a_1+\cfrac{1}{a_2+\cdots}}=[a_0;a_1,a_2,\ldots], \,\,\, a_0 \in \Z, \ a_m \in \N \text{ for } m \geq 1.
\end{equation*}
By $\alpha_m$ we denote the tail of continued fraction $$\alpha_m=  [a_m;a_{m+1},a_{m+2},\ldots].$$
Similarly for $\beta, \gamma \not\in \Q$, let
\begin{align*}
	&\beta = [b_0; b_1,b_2,\ldots] \quad \text{and} \quad \beta_s = [b_s; b_{s+1}, b_{s+2}, \ldots], \\
	&\gamma = [c_0; c_1,c_2,\ldots] \quad \text{and} \quad \gamma_l = [c_l; c_{l+1}, c_{l+2}, \ldots].
\end{align*}
We denote the denominators of convergents of $\alpha, \beta, \gamma$ by $q_m, h_s, r_l$ respectively.
A well-known property of the denominators of convergents is the recursive relation
\begin{equation}
	\label{intro:eq:recur}
	q_{m+1} = a_{m+1}q_m + q_{m-1}.
\end{equation}
Define 
\begin{equation}
	\label{eq:4}
	\xi_m = \psi_\alpha( q_m ) \text{\,\,\,\,\, and \,\,\,\,\,\,} \eta_s = \psi_\beta ( h_s ).
\end{equation}
It is well-known that 
\begin{equation}\label{eq_perron}
	\xi_m  = \frac{1}{q_m\alpha_{m+1}+q_{m-1}}=\frac{1}{q_{m+1}+\frac{q_m}{\alpha_{m+2}}}.
\end{equation}
We will need the following auxiliary results. Lemma \ref{lm:MM} and Lemma \ref{lm:sim_jump} are taken from \cite{Rud22}, and Lemma \ref{lm:3_func} is new and we give a complete proof.
\begin{lemma}[\cite{Rud22} Remark 6]
	\label{lm:MM}
	Let $\alpha, \beta$ be such that $\alpha \pm \beta \not\in \Z$. There exists $T \in \R_+$ such that for all $q_m, h_s > T$ if $$q_{m+1} = h_{s+d}, \ \ h_{s-1} \leq q_{m} < h_s, \ \ \psi_{\alpha}(q_{m+1}-1) < \psi_{\beta}(q_{m+1}-1), $$
	then
	\begin{equation*}
		\psi_{\alpha}(q_m-1) > \psi_{\beta}(q_{m}-1).
	\end{equation*}
	Here $\psi_{\alpha}(q_m - 1)$ is the value of the function before the jump at $q_m$.
\end{lemma}
\import{./}{pic_6.tex}
For an $n$-tuple $\blalpha=(\alpha_1, \dots, \alpha_n)$ of pairwise independent  numbers let
$$\tau(t) = |\{j \in \{1 \dots n \}: \psi_{j} \text{ is discontinuous at } t\}|$$
be the number of functions that have a jump at $t$.
\begin{lemma}[\cite{Rud22} Lemma 7]
	\label{lm:sim_jump}
	There exists $T \in \R_+$ such that for all $t > T$ we have $$\tau(t) \leq \knm(\blalpha).$$
\end{lemma}
\import{./}{pic_1.tex}
\noindent We will also need the following result. 
\begin{lemma}
	\label{lm:3_func}
	Let $\alpha, \beta, \gamma$ be three pairwise independent numbers with corresponding functions $\psi_{\alpha}, \psi_{\beta}, \psi_{\gamma}$. If
	\begin{align*}
		q_m &= h_s, & q_{m+1} &= r_l, & r_{l+1} & = h_{s+1}, \\
		q_{m+2} &= h_{s+2}, & q_{m+3} &= r_{l+2},  & r_{l+3} &= h_{s+3},
	\end{align*}
	then
	\begin{equation}
		\label{eq:5}
		\psi_{\beta}(t) > \psi_{\alpha}(t), \ \forall t \, \in [h_{s+1}, h_{s+2}).
	\end{equation}    
\end{lemma}
\begin{proof}
	Using (\ref{eq:4}) we rewrite (\ref{eq:5}) as
	\begin{equation*}
		\eta_{s+1} > \xi_{m+1}.
	\end{equation*}
	From the growth of the denominators of convergences we immediately establish the inequalities
	\begin{equation*}
		q_m = h_s < q_{m+1} = r_l < r_{l+1} = h_{s+1} < h_{s+2} = q_{m+2} < q_{m+3} = r_{l+2} < r_{l+3} = h_{s+3}.
	\end{equation*} 
	We will prove it by contradiction. Suppose that $\xi_{m+1} > \eta_{s+1}$. From Perron's formula (\ref{eq_perron}) we have
	\begin{equation*}
		\frac{1}{q_{m+1}\alpha_{m+2}+q_m} > \frac{1}{h_{s+1}\beta_{k+2}+h_s},
	\end{equation*}
	and since $q_m = h_s$ we get
	\begin{equation*}
		h_{s+1} \beta_{k+2} > q_{m+1} \alpha_{m+2},
	\end{equation*}
	or by definitions of $\alpha_{m+2}$ and $\beta_{k+2}$ we can rewrite it as
	\begin{equation}
		\label{lm_3func_eq1}
		h_{s+1} \left(b_{s+2} + \frac{1}{\beta_{s+3}}\right) > q_{m+1} \left(a_{m+2} + \frac{1}{\alpha_{m+3}}\right).
	\end{equation}
	Now from the two equalities $h_{s+2} = q_{m+2}$ and $h_s = q_m$ by the recursive formula for denominators (\ref{intro:eq:recur}) from
	\begin{gather*}
		q_{m+1} a_{m+2} + q_m = h_{s+1} b_{s+2} + h_s
	\end{gather*}
	we get
	\begin{equation*}
		q_{m+1} a_{m+2} = h_{s+1} b_{s+2}.
	\end{equation*}
	We substitute this result to (\ref{lm_3func_eq1}) and get
	\begin{equation}
		\label{lm_3func_eq3}
		h_{s+1} \alpha_{m+3} > q_{m+1} \beta_{s+3}.
	\end{equation}
	Because $r_{l+3} = h_{s+3}$ and $r_{l+1} = h_{s+1}$, again by (\ref{intro:eq:recur}) 
	\begin{equation*}
		c_{l+3} r_{l+2} + r_{l+1} = b_{s+3} h_{s+2} + h_s
	\end{equation*}
	we get
	\begin{equation}
		\label{lm_3func_eq2}
		c_{l+3} r_{l+2} = b_{s+3} h_{s+2}.
	\end{equation}
	Since $h_{s+2} = q_{m+2}$ and $r_{l+2} = q_{m+3} = a_{m+3} q_{m+2} + q_{m+1}$ equation (\ref{lm_3func_eq2}) becomes
	\begin{gather*}
		c_{l+3}(a_{m+3} q_{m+2} + q_{m+1}) = b_{s+3} q_{m+2},
	\end{gather*}
	or
	\begin{equation*}
		c_{l+3} q_{m+1} = q_{m+2}(b_{s+3}-c_{l+3}a_{m+3}).
	\end{equation*}
	Because two consecutive denominators are relatively prime, there exists $d \in \N$ such that 
	\begin{align*}
		&c_{l+3} = d q_{m+2}, \\
		&b_{s+3}-c_{l+3}a_{m+3} = dq_{m+1}
	\end{align*}
	(note that because $c_{l+3}$ and $q_{m+2}$  are positive, so is $d$). Since $\beta_{s+3} > b_{s+3}$ we can substitute this to (\ref{lm_3func_eq3}) and get
	\begin{equation*}
		h_{s+1} \alpha_{m+3} >q_{m+1} (dq_{m+1} + c_{l+3}a_{m+3}) =  dq_{m+1} (q_{m+1} + q_{m+2}a_{m+3}).
	\end{equation*}
	For the left-hand side of this equation we use inequalities $q_{m+2} = h_{s+2} > h_{s+1}$ and $a_{m+3} + 1 > \alpha_{m+3}$ to obtain
	\begin{gather*}
		q_{m+2}(a_{m+3}+1) > dq_{m+1} (q_{m+1} + q_{m+2}a_{m+3}), 
	\end{gather*}
	or
	\begin{equation*}
		q_{m+2}(a_{m+3}+1 - dq_{m+1}a_{m+3}) > dq_{m+1}^2 > 0.
	\end{equation*}
	However, $a_{m+3}+1 < dq_{m+1}a_{m+3}$ for all $m \geq 1$, so the left-hand side is negative, which is a contradiction. This proves Lemma \ref{lm:3_func}.
\end{proof}
\begin{remark}
	Lemma \ref{lm:3_func} does not hold for only two functions $\psi_\alpha$ and $\psi_\beta$.
\end{remark}
\subsection{The permutation $\pi$}
\begin{lemma}
	\label{lm:order}
	The order of permutation $\pi$ is $k$.
\end{lemma}
\begin{proof}
	From Section \ref{subsec:pi} recall the definition of the permutation $\pi$ in terms of a $k \times k$ diagram (see Figure \ref{fig:pi}). The permutation $\pi$ acts as follows (see Figure \ref{fig:pi:sub1}):
	\begin{itemize}
		\item the first row is shifted cyclically to the left, where the first element moves to the last place in the row;
		\item the first column without the top element is placed horizontally under the first row in reverse order;
		\item the remaining elements are shifted down one row together.
	\end{itemize}
	\import{./}{pic_11.tex}	 
	It is clear that the permutation restricted to the first row is a cyclic shift of order $k$. We now consider the remaining sub-diagram $(k-1) \times (k-1)$, which is obtained from the initial diagram by removing the first row (see Figure \ref{fig:pi:sub2}). The $d$th element $A$ from the first column (counting from the top) is moved to the $d$th position $B$ in the first row (counting from the right). All the other elements are shifted along the diagonal direction (down and to the left). For example, the element $C$ moves to $D$. 
	
	There are two cases: $k$ is odd and $k$ is even. 
	
	In the case $k$ is odd, the orbit of any element consists of a pair of diagonals. The first diagonal is the one to which the element belongs. Let $d$ denote the position of the element from the first column within this diagonal (counting from the top). Then the second diagonal is the one to which $(k-d)$th element from the first column belongs (counting from the top).  
	
	In the case $k$ is even, for $d = k / 2$ two diagonals coincide. So, the orbit of any element in the middle diagonal consists of only this middle diagonal, and the order of such elements is $k/2$. The orbit of any other element  again consists of two diagonals as in the first case.
	
	The sum of the lengths of two diagonals belonging to one orbit is $k$, so the order of any element (except for elements in the middle diagonal when $k$ is even) is $k$. 
	This completes the proof.
\end{proof}
\begin{remark}
	From the proof of Lemma \ref{lm:order} it follows that when $k$ is odd, the permutation $\pi$ can be decomposed into a product of $(k+1)/2$ cycles, and when $k$ is even it can be decomposed into a product of $(k+2)/2$ cycles.
\end{remark}
	The proof of the following proposition is a direct computation, and is left to the reader.
\begin{proposition}
	\label{rm:pi}
	 We put the following vector in the diagram, 
	\begin{multline*}
		\blw = \left(u_{k,k}, \, u_{k-1,k},\, u_{k-2,k},\, \dots,\, u_{1,k}, \, u_{1,1},\, u_{1,k-1},\, u_{1,k-2}, \dots,\right.  \\ \left.\dots, \, u_{1,2}, \, u_{2,2}, \, u_{2,k-1}, \dots, \, u_{2,3}, \, u_{k-2,k-2},\, u_{k-2,k-1},\, u_{k-1,k-1}\right),
	\end{multline*}
	and get
	\begin{equation*}
		\pi(\blw) = (u_{1,1}, \, u_{1,k},\, u_{1,k-1},\, \dots,\, u_{1,2}, \, u_{2,2},\, u_{2,k}, \dots, \, u_{2,3}, u_{3,3}, \dots, \,u_{k-1,k-1},\, u_{k-1,k},\, u_{k,k}).
	\end{equation*}
\end{proposition}

\subsection{Projection of the order vector \texorpdfstring{$\boldsymbol{v_{\alpha}}(t)$}{v}}
For the collection of functions $\psi_1, \dots, \psi_n$ which corresponds to the $n$-tuple $\a = (\alpha_1, \dots, \alpha_n)$ we consider a certain subcollection of functions $\psi_{i_1}, \dots, \psi_{i_m}$ which corresponds to the $m$-tuple 
\begin{equation}
	\label{eq:39}
	\bbeta = (\alpha_{i_1}, \dots, \alpha_{i_m}).
\end{equation}
Let $\mathcal{V}_{\a}$
be the set of values of $\blv_{\a}(t)$ which occur infinitely many times:
\begin{equation*}
	\mathcal{V}_{\a} = \{\blv = (v_1, 
	\dots, v_n) \in S_n: \  \forall \, t > 0 \ \exists \, t' > t \ \blv = \blv_{\a}(t')\}.
\end{equation*}
Let $I = (i_1, \dots, i_m)$ denote the indices of the subcollection $\bbeta = (\alpha_{i_1}, \dots, \alpha_{i_m})$. By $\Sym(I)$ we denote the permutation group of $I$. We now define a projection operator
\begin{align*}
	\pr_{\bbeta}: \mathcal{V}_{\a} \rightarrow \Sym(I), \quad
	\pr_{\bbeta}\bigl(\blv=(v_1, \dots, v_n)\bigr) = (v_{r_1}, \dots, v_{r_m}),
\end{align*}
where $r_1 < r_2 < \dots < r_{m-1} < r_m$ and for all $1 \leq s \leq m$ there exists $1 \leq j \leq m$ such that $v_{r_s} = i_j$. That is, the image of $\pr_{\bbeta}(\blv)$ is a vector $\blu \in \Sym(I)$ which consist of indices $(i_1, \dots, i_m)$ in the same order as they are in $\blv$.
The preimage of $\blu \in \Sym(I)$ is defined as
\begin{equation}
	\label{eq:41}
	\pr^{-1}(\blu) = \{\blv \in \mathcal{V}_{\a}: \ \pr_{\bbeta}(\blv) = \blu\}
\end{equation}
The operator $\pr_{\bbeta}$ is not injective, so the preimage may be a set of several vectors.
For convenience we will slightly abuse the notation of an order vector. We define an order vector on the subcollection $\bbeta$ as
\begin{equation}
	\label{eq:40}
	\blv_{\bbeta}(t) = \pr_{\bbeta}(\blv_{\a}(t)).
\end{equation}
We also note that clearly for different $\blu_1, \blu_2 \in \Sym(I)$
\begin{equation*}
	\pr_{\bbeta}^{-1}(\blu_1) \cap \pr_{\bbeta}^{-1}(\blu_2) = \emptyset. 
\end{equation*}
To clarify the definitions, we consider an example. 
\begin{example}\label{ex:1} Let $\a = (\alpha_1, \alpha_2, \alpha_3, \alpha_4)$ and $\bbeta = (\alpha_3, \alpha_4)$. Thus,
	\begin{equation*}
		I = (3, 4), \quad \Sym(I) = \bigl\{(3, 4), (4, 3) \bigr\}.
	\end{equation*}
	Let 
	\begin{equation*}
		\mathcal{V}_{\a} = \bigl\{ (1,2,3,4), \ (3,2,4,1), \ (4, 1, 2, 3) \bigr\}.
	\end{equation*} The projection operator acts as 
	\begin{align*}
		&\pr_{\bbeta}\bigl((1,2,3,4)\bigr) = (3,4) ,&\ &\pr_{\bbeta}\bigl((3, 2, 4, 1)\bigr) = (3,4), &\
		&\pr_{\bbeta}\bigl((4, 1, 2, 3)\bigr) = (4,3).
	\end{align*}
	The preimages are
	\begin{equation*}
		\pr_{\bbeta}^{-1}\bigl((3,4)\bigr) = \bigl\{ (1,2,3,4), \ (3, 2, 4, 1) \bigr\}, \quad \pr_{\bbeta}^{-1}\bigl((4,3)\bigr) = \bigl\{ (4, 1, 2, 3) \bigr\}.
	\end{equation*}
\end{example}
\import{./}{pic_9.tex}

	In the proof of Theorem \ref{th:main} \ref{it2} we will use the argument about existence of several vectors in the preimage of projected vector $\blv_{\bbeta}(t_i)$ for a subcollection of continuous functions at  $t_i$. Here we describe this idea on Example \ref{ex:1}. Functions $\psi_{3}$ and $\psi_4$  in Figure \ref{fig:pr} are both continuous at $t_{i+1}$, so their relative vector order $\blv_{\bbeta}(t)$ does not change value at $t_{i+1}$, that is
	\begin{equation*}
		\blv_{\bbeta}(t_i) = \pr_{\bbeta}(\blv_{\a}(t_{i})) = \pr_{\bbeta}(\blv_{\a}(t_{i+1})) = \blv_{\bbeta}(t_{i+1}) = (3,4).
	\end{equation*}
	However, $t_i$ is a moment of value change for $\blv_{\a}(t)$,
	\begin{equation*}
		\blv_{\a}(t_i) \neq \blv_{\a}(t_{i+1}).
	\end{equation*}
	So, the preimage consists of at least two vectors,
	\begin{align*}
		&\{ \blv_{\a}(t_{i}), \blv_{\a}(t_{i+1})\} \subset \pr_{\bbeta}^{-1}\bigl(\blv_{\bbeta}(t_i)\bigr), \\
		&\bigl\{ (1,2,3,4), \ (3, 2, 4, 1)\bigr\} \subset \pr_{\bbeta}^{-1}\bigl((3,4)\bigr).
	\end{align*}

\section{A strengthening of Theorem \ref{th:uniqness}} \label{sec:4}
For $t_0 \in \R_{+}$ we renumerate the functions $\left\{\psi_i\right\}_{i=1}^n$ with the indices $(i,j)$ as $\left\{\psi_{i,j}\right\}_{1\le i \le j \le k}$ such that at $t_0$ they are ordered as
\begin{multline}
	\label{eq:2}
	\psi_{1,1}(t_0) > \psi_{1,k}(t_0) > \psi_{1,k-1}(t_0) > \dots > \psi_{1,2}(t_0) > \\ >  \psi_{2,2}(t_0) > \psi_{2,k}(t_0) > \psi_{2,k-1}(t_0) > \dots > \psi_{2,3}(t_0) > \\ > \psi_{3,3}(t_0)  > \dots > \psi_{k-1}(t_0) > \psi_{k-1, k}(t_0) > \psi_{k,k}(t_0),
\end{multline}
in accordance with (\ref{int:eq10}).
With this notation in hand we now state a strengthening of Theorem \ref{th:uniqness}.
\begin{theorem}
	\label{th:main}
	Let $\a = (\alpha_1, \dots, \alpha_n)$ be a $\frac{k(k+1)}{2}$-tuple of pairwise  independent real numbers with $\knm(\a) = k$.  Then there exists $T \in \R_{+}$ such that for any $t_0 \geq T$ for the sequence $\{t_i\}_{i=1}^\infty = \{t_i(\a, t_0)\}_{i=1}^\infty$ of moments of permutation changes defined in (\ref{def:ti}) and for the enumeration (\ref{eq:2}) the following holds.
	\begin{enumerate}[label=(\roman*)]	
		\item \label{it1} For all $i \geq 0$
		\begin{equation*}
			\tau(t_i) = k.
		\end{equation*} 
		That is, at each $t_i,\ i \geq 1$ exactly $k$ functions have a jump.
		\item \label{it2} For all $\ 0 \leq i \leq k-1, \, s \geq 0$ 
		\begin{equation*}
			\blv_{\a}(t_{i+sk}) = \blv_{\a}(t_i).
		\end{equation*}
		This means that the sequence $\{\blv_{\a}(t_i)\}_{i=0}^\infty$  is periodic with period $k$.
		\item \label{it3} Each function $\psi_{i,i}, \ 1 \leq i \leq k$ has jumps at $\{t_{i+sk}\}_{s=0}^\infty$ and is continuous at all other $\{t_j\}_{j=1}^\infty$.
		\item \label{it4} Each function $\psi_{i,j}, \ 1 \leq i \leq k, \ i < j \leq k$ has jumps at $\{t_{i+sk}\}_{s=0}^\infty$ and $\{t_{j+sk}\}_{s=0}^\infty$, and is continuous everywhere between these jumps.
		\item \label{it5} For all $ i \geq 1, \, s \in \N_0,$ the
		$k$ largest functions at $t_{i-1+ks}$ have a jump at $t_{i+ks}$. The function $\psi_{i, i}$ is the largest function at $t_{i-1+ks}$.
		\item \label{it6} For all $i \geq 0$
		\begin{equation*}
			\blv_{\a}(t_{i+1}) = \pi(\blv_{\a}(t_i)).
		\end{equation*}
	\end{enumerate}
\end{theorem}
Theorem \ref{th:uniqness} follows from Theorem \ref{th:main} \ref{it6}. Figures \ref{fig:ex3} and \ref{fig2} show the cases $k=2,\,n=3$ and  $k=3, \, n=6$ respectively.

	Here we note the reasoning for the enumeration (\ref{int:eq10}). The index $(i,i)$ of the function $\psi_{i,i}$ corresponds to the jumps the function $\psi_{i,i}$ has at $\{t_{i+sk}\}_{s=0}^\infty$ by \ref{it3}. The index $(i,j)$ of the function $\psi_{i,j}$ corresponds to the jumps of the function $\psi_{i,j}$ at $\{t_{i+sk}\}_{s=0}^\infty$ and $\{t_{j+sk}\}_{s=0}^\infty$ by \ref{it4}. 

	Note that we do not claim that functions $\psi_{i,i}$ are continuous between the jumps $\{t_{i+sk}\}_{s=0}^\infty$, only that they are continuous at all other points $\{t_i\}_{i=0}^\infty$.

\section{Proof of Theorem \ref{th:main}} \label{ch_2:uniqness}

\subsection{Outline of the proof} Item \ref{it1} follows from Theorem \ref{th:upper_bound} and Lemma \ref{lm:sim_jump}. Items \ref{it2} -- \ref{it6} are proved by induction, with \ref{it3} -- \ref{it6}  treated collectively.

The main idea of the proof of  items \ref{it2}--\ref{it6} is as follows. Let $\a = (\alpha_1, \dots, \alpha_n)$ be  $\frac{k(k+1)}{2}$-tuple of pairwise  independent real numbers with $\knm(\a) = k$. As the induction hypothesis, we assume that for any collection $\bbeta$ of $\frac{k(k-1)}{2}$ independent numbers with $\knm(\bbeta) = k-1$ the statement of the theorem holds. We then consider a certain subcollection of functions $\psi_{i_1}, \dots, \psi_{i_m}$ that are continuous at some $t_i$ and correspond to an $m$-tuple of real numbers
$\bbeta_i = (\alpha_{i_1}, \dots, \alpha_{i_m})$. By \ref{it1}, the size of the subcollection $\bbeta_i$ is $\frac{k(k-1)}{2}$. 
Because of the continuity of these functions, there are successive equal vectors $\blv_{\bbeta_i}(t_{i-1}) = \blv_{\bbeta_i}(t_i).$ Thus, the number of distinct values of $\blv_{\bbeta_i}(t)$ is smaller than the number of distinct values of $\blv_{\a}(t)$ at least by one, that is $$\knm(\bbeta_i) \leq \knm(\a) - 1.$$
We then apply the induction hypothesis to the subcollection $\bbeta_i$.

In the proof, depending on convenience, we use both enumerations of functions: the natural enumeration $\psi_1, \dots, \psi_n$ and enumeration (\ref{eq:2}).

For a given $n$-tuple $\a = (\alpha_1, \dots, \alpha_n)$ of independent real numbers we first set $T$ as in definition (\ref{def:T}), so that $\blv_{\a}(t)$ is correctly defined for all $t \geq T$. Since we are proving the existence of $T$, at some point in the proof, if needed, we can set $T$ larger. We fix an arbitrary $t_0 \geq T$ and construct the sequence $\{t_i\}_{i=1}^{\infty} =\{t_i(\a, t_0)\}_{i=1}^{\infty}$ by (\ref{def:ti}).
\subsection{Proof of \ref{it1}}
We prove the statement for $t_1$; the same argument applies to all $t_i$. By assumption, there are $k$ distinct values of vector $\blv_{\a}(t)$ that appear infinitely many times as $t \rightarrow \infty$. We denote them by $\blv_1, \dots, \blv_k$, with $\blv_1 = \blv_{\a}(t_0), \ \blv_2 = \blv_{\a}(t_1)$. Let $\psi_{i_1}, \dots, \psi_{i_m}$ be the functions that are continuous at $t_1$, and $\bbeta = (\alpha_{i_1}, \dots, \alpha_{i_m})$ be the collection of numbers corresponding to these functions. Since these functions are continuous at $t_1$, their vector $\blv_{\bbeta}(t)$ does not change at $t_1$, that is
\begin{equation}
	\label{eq14}
	\blv_{\bbeta}(t_0) = \blv_{\bbeta}(t_1).
\end{equation}
Since $\blv_{\a}(t)$ has $k$ distinct values, the vector $\blv_{\bbeta}(t)$ has at most $k$ distinct values. Furthermore, since two of them are equal by (\ref{eq14}), we get $\knm(\bbeta) \leq k-1$. From this and Theorem \ref{th:upper_bound} it follows that
\begin{equation*}
	m \leq \frac{(k-1)k}{2}.
\end{equation*}
By Lemma \ref{lm:sim_jump} no more than $k$ functions have a jump at $t_1$, thus
\begin{equation*}
	m \geq \frac{k(k+1)}{2} - k =\frac{(k-1)k}{2}.
\end{equation*}
Hence, we have $$m = \frac{k(k-1)}{2},$$ and so $$\tau(t_1) = \frac{k(k+1)}{2} - m = k.$$
This completes the proof of \ref{it1}.

\subsection{Proof of \ref{it2}.} 
We prove that the sequence $\{\blv_{\a}(t_i)\}_{i=0}^\infty$  is periodic and $k$ is the period. We write the sequence of vectors $\{\blv_{\a}(t_i)\}_{i=0}^\infty$ as an infinite word 
\begin{equation*}
	W = \blv_{\a}(t_0) \blv_{\a}(t_1)\blv_{\a}(t_2)\blv_{\a}(t_3)\blv_{\a}(t_4)\blv_{\a}(t_5) \dots.
\end{equation*}
In this notation, we want to prove that the first $k$ values $\blv_{\a}(t_i), \ 0 \leq i \leq k-1,$ are all different, and that our word is periodic, that is
\begin{equation*}
	W = \overline{\blv_{\a}(t_0) \blv_{\a}(t_1)\blv_{\a}(t_2)\blv_{\a}(t_3) \dots \blv_{\a}(t_{k-2}) \blv_{\a}(t_{k-1}) \vphantom{\dfrac{a}{b}}}.
\end{equation*}
As was noted at the end of the outline of the proof, we are proving the existence of $T$, so it is actually enough to prove that the word $W$ is eventually periodic, 
\begin{equation*}
	W = \blv_{\a}(t_0) \blv_{\a}(t_1)\blv_{\a}(t_2)\blv_{\a}(t_3) \dots \blv_{\a}(t_{s-1}) \, \overline{\blv_{\a}(t_{s}) \dots \blv_{\a}(t_{s+k-1}) \vphantom{\dfrac{a}{b}}}.
\end{equation*}
and set $T = t_{s}$.

In the case $k = 2, \, n = 3$ there are only two values of the vector $\blv_{\a}(t)$, so there is nothing to prove.

For $k \geq 3$ we prove \ref{it2} by induction, starting with $k=3$. 

We note that the proof of the base case $k = 3$ requires only five functions, and not six.
%Step of induction does not work unless $k \geq 4$ (see Remark \ref{rm:3}).

\textbf{The base case}. Let $k=3$ and $n=6$. Let $\a=(\alpha_1, \dots, \alpha_6)$ be six pairwise independent numbers with $\knm(\a) = 3$ and let $\blv_1, \blv_2, \blv_3$ denote three distinct values of the vector $\blv_{\a}(t)$. We want to show that the word $W$ is eventually periodic, 
\begin{equation*}
	W = \blv_i \blv_j \dots \blv_l \, \overline{\blv_1\blv_2\blv_3 \vphantom{T}}.
\end{equation*}
\import{./}{pic_3.tex}
Since all three values appear infinitely many times, to prove this, it is enough to show that $W$ does not contain infinitely many subwords of the form $abac$, where $a,b,c$ is an arbitrary combination of different letters $\blv_1, \blv_2, \blv_3$. 

 We assume the contrary, namely, that  subword  $abac$ appears infinitely many times. Let $a = \blv_1, b = \blv_2, c = \blv_3$.

Suppose that $\blv_1\blv_2\blv_1\blv_3$ appears at  $t_i$, that is, $$\blv_{\a}(t_i) = \blv_1,\ 
\blv_{\a}(t_{i+1}) = \blv_2, \ \blv_{\a}(t_{i+2}) = \blv_1, \ \blv_{\a}(t_{i+3}) = \blv_3, \ \blv_{\a}(t_{i+4}) = \blv_{1/2},$$ where $\blv_{1/2}$ is the notation for either $\blv_1$ or $\blv_2$ (see Figure \ref{fig:period}).

By \ref{it1}, three functions are continuous at $t_{i+1}$.
Without loss of generality we assume that these functions are $\psi_1, \psi_2, \psi_3$, and for all $t \in [\,t_i, \,t_{i+3})$ (since they are continuous at $t_{i+1}$ and at $t_{i+2}$ the vector is again $\blv_1$), we have
\begin{equation*}
	\psi_{1}(t) > \psi_{2}(t) > \psi_3(t).
\end{equation*} 
So, for the subcollection $(\alpha_1, \alpha_2, \alpha_3)$ we have
\begin{equation*}
	\blv_{(\alpha_1, \alpha_2, \alpha_3)}(t_{i}) = \blv_{(\alpha_1, \alpha_2, \alpha_3)}(t_{i+1}),
\end{equation*}
and $\knm\bigl((\alpha_1, \alpha_2, \alpha_3)\bigr)=2$.
By Kan-Moshchevitin Theorem \ref{th:KM}, any two functions interchange infinitely many times, so the order of $\psi_1, \, \psi_2, \, \psi_3$ at $t_{i+3}$ should be reversed:
\begin{equation*}
	\psi_{3}(t_{i+3}) > \psi_{2}(t_{i+3}) > \psi_1(t_{i+3}).
\end{equation*} 
Hence, $\psi_1$ and $\psi_2$ have a jump at $t_{i+3}$. Then, by Lemma \ref{lm:MM} $\psi_2$ is continuous for all $t \in [t_{i}, t_{i+3})$. By assumption, $\blv_{\alpha}(t_{i+4})$  is either $\blv_1$ or $\blv_2$, so 
\begin{equation*}
	\psi_{1}(t_{i+4}) > \psi_{2}(t_{i+4}) > \psi_3(t_{i+4}),
\end{equation*} 
and hence $\psi_2$ and $\psi_3$ have a jump at $t_{i+4}$.

By \ref{it1}, three functions are continuous at $t_{i+3}$. We keep in mind that the collection of functions $\psi_1, \psi_2, \psi_3$ has only two values of $\blv_{\alpha}(t)$, that is, $\knm\bigl((\alpha_1, \alpha_2, \alpha_3)\bigr) = 2$. So, by Lemma \ref{lm:sim_jump}, only two of these functions can have a simultaneous jump. Since $\psi_1$ and $\psi_2$ have jump at $t_{i+3}$, function $\psi_3$ is continuous at $t_{i+3}$. Let $\psi_{4}$ and $\psi_{5}$ denote the other two functions that are continuous at $t_{i+3}$. We can assume that $\psi_4(t_i) > \psi_5(t_i)$.

Functions $\psi_3, \psi_4, \psi_5$ are continuous at $t_{i+3}$, so 
\begin{equation*}
	\blv_{(\alpha_3, \alpha_4, \alpha_5)}(t_{i+2}) = \blv_{(\alpha_3, \alpha_4, \alpha_5)}(t_{i+3}).
\end{equation*}
Then by Kan-Moshchevitin Theorem \ref{th:KM}, the vector $\blv_{(\alpha_3, \alpha_4, \alpha_5)}(t_{i})$  is reversed with respect to the vector $\blv_{(\alpha_3, \alpha_4, \alpha_5)}({t_{i+1}})$. From this and because $\psi_3$ is continuous at $t_{i+1}$, functions $\psi_{4}$ and $\psi_5$ are greater than $\psi_{3}$ at $t_i$, and have a jump under $\psi_3$ at $t_{i+1}$, that is, 
\begin{equation*}
	\psi_{4}(t_i) > \psi_{5}(t_i) > \psi_3(t_i),
\end{equation*}
and 
\begin{equation}
	\label{eq:6}
	\psi_{3}(t_{i+1}) > \psi_{5}(t_{i+1}) > \psi_4(t_{i+1}).
\end{equation}
From (\ref{eq:6}) and $\psi_2(t_{i+1}) > \psi_3(t_{i+1})$, we get
\begin{equation}
	\label{eq:3}
	\psi_{2}(t_{i+1}) > \psi_4(t_{i+1}) > \psi_{5}(t_{i+1}).
\end{equation}
Because $\psi_2$ is continuous at $t_{i+1}$ and $t_{i+2}$, and $\blv_{\alpha}(t_{i+2}) = \blv_1$, function $\psi_4$ cannot have a jump under $\psi_{2}$ at $t_{i+1}$. 
That is, we cannot have $$\psi_4(t_i) > \psi_2(t_i) \text{ and } \psi_2(t_{i+1}) > \psi_4(t_{i+1}).$$ So, by (\ref{eq:3}), we get
\begin{equation*}
	\psi_{2}(t_{i}) > \psi_4(t_{i}) > \psi_{5}(t_{i}).
\end{equation*}
Then, by Kan-Moshchevitin Theorem \ref{th:KM}, $\psi_2$ has a jump under $\psi_4$ and $\psi_5$ at $t_{i+3}$, that is,
\begin{equation*}
	\psi_4(t_{i+3}) > \psi_5(t_{i+3}) > \psi_2(t_{i+3}).
\end{equation*}
Since $\blv_{\alpha}(t_{i+4})$ 
is either $\blv_1$ or $\blv_2$, at $t_{i+4}$ we have
\begin{equation*}
	\psi_{2}(t_{i+4}) > \psi_4(t_{i+4}) > \psi_{5}(t_{i+4}),
\end{equation*}
or
\begin{equation*}
	\psi_{2}(t_{i+4}) > \psi_5(t_{i+4}) > \psi_{4}(t_{i+4}).
\end{equation*}
In both cases $\psi_4$ and $\psi_5$ should have a jump at $t_{i+4}$ under $\psi_2$.

We have shown that four functions $\psi_2, \psi_3, \psi_4, \psi_5$ have a jump at $t_{i+4}$.
Because $\knm(\alpha) = 3$, this is a contradiction with Lemma \ref{lm:sim_jump}. This completes the proof of the base of case.

\textbf{Step of induction.}  Recall from definition (\ref{def:ti}) that the sequence $\{t_i(\a, t_0)\}_{i=0}^\infty$ depends on $\a$. In this step of the proof, we will consider two subcollections of numbers $\bbeta_1$ and  $\bbeta_2$ of the form (\ref{eq:39}), for which the sequences $\{t_i(\bbeta_j, t_0)\}_{i=0}^\infty, \ j =1, 2$ may be different from $\{t_i(\a, t_0)\}_{i=0}^\infty.$ 

Let $k \geq 4$. As the induction hypothesis we assume that for any $\frac{k(k-1)}{2}$-tuple $\bbeta$ of pairwise independent numbers with $\knm(\bbeta) = k-1$ there exist $T = T(\bbeta) \in \R_+$ such that for any $t_0 \geq T$ vectors $\blv_{\bbeta}(t_i(\bbeta)), \ 0 \leq i \leq k-2$ are all different and the sequence $\{\blv_{\bbeta}(t_i(\bbeta, t_0))\}_{i=0}^\infty$  is periodic with period $k-1$. 

Let $\a$ be a $\frac{k(k+1)}{2}$-tuple of pairwise independent numbers with $\knm(\a) = k$. Using the induction hypothesis we want to find $T = T(\a)$ such that for all $t_0 \geq T$ the vectors $\blv_{\a}(t_i(\a, t_0)), \ 0 \leq i \leq k-1$ are all different and  the sequence $\{\blv_{\a}(t_i(\a, t_0))\}_{i=0}^\infty$ is periodic with period $k$. 

As before, we first set $T = T(\a)$ as in definition (\ref{def:T}), so that $\blv_{\a}(t)$ is correctly defined for all $t \geq T$. We will find two subcollections $\bbeta_1$ and $\bbeta_2$ of $\a$, to which we will apply the induction hypothesis, obtaining $T(\bbeta_1)$ and $T(\bbeta_2)$. Since we are proving the existence of $T$, we can set $T = \max\bigl(T(\bbeta_1), T(\bbeta_2)\bigr)$, and so for all $t_0 \geq T$ sequences $\{\blv_{\bbeta_1}(t_i(\bbeta_1, t_0))\}_{i=0}^\infty$ and $\{\blv_{\bbeta_2}(t_i(\bbeta_2, t_0))\}_{i=0}^\infty$ are periodic.

Let $t_0 \geq T$ and $\{t_i\}_{i=1}^{\infty} =\{t_i(\a, t_0)\}_{i=1}^{\infty}$ as in (\ref{def:ti}).
By \ref{it1} applied to $\a$ there are $m = \frac{k(k-1)}{2}$ functions $\psi_{i_1}, \dots, \psi_{i_m}$ which are continuous at $t_{1}(\a)$. We denote the  subcollection of numbers corresponding to  functions $\psi_{i_1}, \dots, \psi_{i_m}$ by $\bbeta_1 = (\alpha_{i_1}, \dots, \alpha_{i_m})$. Since for all functions we have $\knm(\a) = k$, for a subcollection of functions we also have $\knm(\bbeta_1) \leq k$. Recall definition (\ref{eq:40}) of a vector for a subcollection. Since functions $\psi_{i_1}, \dots, \psi_{i_m}$ are continuous at $t_{1}(\a)$, they do not change their relative order, thus we have two successive equal vectors
\begin{equation*}
	\blv_{\bbeta_1}(t_0) = \blv_{\bbeta_1}(t_{1}(\a)).
\end{equation*} 
At the same time, for the entire collection of functions, by definition (\ref{def:ti}), we have $\blv_{\a}(t_{0}) \neq \blv_{\a}(t_{1}(\a))$. So,  $$\knm(\bbeta_1) \leq \knm(\a) - 1 = k - 1.$$ By Theorem \ref{th:upper_bound}, since $m = \frac{k(k-1)}{2}$, we necessarily have 
\begin{equation*}
	\knm(\bbeta_1) = k-1.
\end{equation*}
By induction hypothesis applied to $\bbeta_1$ there exists $T(\bbeta_1)$ such that for all $t'_0 \geq T(\bbeta_1)$ we have a $(k-1)$-periodic sequence $\{\blv_{\bbeta_1}(t_i(\bbeta_1, t'_0))\}_{i=0}^\infty$, where $\blv_{\bbeta_1}(t_{i}(\bbeta_1, t'_0)),$ $ 0 \leq i \leq k-2$ are all different. 
If $T(\bbeta_1) \geq T$, we set $T = T(\bbeta_1)$ and fix a new arbitrary $t_0 \geq T$.
We denote $$\blu_{i+1} = \blv_{\bbeta_1}(t_i(\bbeta_1, t_0)), \ 0 \leq i \leq k-2, $$ and write this $(k-1)$-periodic sequence as a periodic word 
\begin{equation}
	\label{eq:20}
	W_{\bbeta_1} = \overline{ \blu_1 \blu_2  \dots \blu_{k-2} \blu_{k-1}\vphantom{T}}.
\end{equation}  
We consider the preimage of $\blu_1$ under $\pr_{\bbeta_1}$ (recall the definition of the preimage of the projection operator (\ref{eq:41}) and Example \ref{ex:1}).
Because $\blu_1 = \blv_{\bbeta_1}(t_0) = \blv_{\bbeta_1}(t_{1}(\a))$ and $\blv_{\a}(t_0) \neq \blv_{\a}(t_1)$, 
this preimage consists of at least two vectors $\blv_{\a}(t_0)$ and $\blv_{\a}(t_1)$, that is,
\begin{equation}
	\label{eq:23}
	\Bigl\{\blv_{\a}(t_0), \blv_{\a}\bigl(t_1(\a)\bigr)\Bigr\} \subset \pr^{-1}_{\bbeta_1}(\blu_1).
\end{equation}
Different values of the vector $\blv_{\bbeta_1}(t)$ on the subcollection $\bbeta_1$ correspond to different values of the vector $\blv_{\a}(t)$ on the entire collection $\a$, so the preimages $\pr^{-1}(\blu_i), \ 0 \leq i \leq k-2$ are pairwise disjoint sets. To summarize, we have
\begin{align*}
	&\bigl|\pr^{-1}_{\bbeta_1}(\blu_1)\bigr| \geq 2, \\
	&\bigl|\pr^{-1}_{\bbeta_1}(\blu_i)\bigr| \geq 1, \ 2 \leq i \leq k-1, \\
	&\pr^{-1}_{\bbeta_1}(\blu_i) \cap \pr^{-1}_{\bbeta_1}(\blu_j) = \emptyset, \ i \neq j,
\end{align*}
where $\bigl|\pr^{-1}_{\bbeta_1}(\blu_i)\bigr|$ denotes the number of elements in $\pr^{-1}_{\bbeta_1}(\blu_i)$.
there are exactly $k$ distinct values of $\blv_{\a}(t)$, we get
\begin{align}
	\label{eq:31}
	&|\pr^{-1}_{\bbeta_1}(\blu_1)| = 2, \\
	\label{eq:32}
	&|\pr^{-1}_{\bbeta_1}(\blu_i)| = 1, \ 2 \leq i \leq k-1.
\end{align}
We denote $k$ distinct values of $\blv_{\a}(t)$ by
\begin{align}
	\blv_1 &= \blv_{\a}(t_0), \ \blv_2 = \blv_{\a}\bigl(t_1(\a)\bigr),  \label{eq:21} \\
	\blv_i &= \pr^{-1}_{\bbeta_1}(\blu_{i-1}), \ 3 \leq i \leq k. \label{eq:22}
\end{align}
We replace $\blu_i$ in (\ref{eq:20}) by preimages (\ref{eq:21}), (\ref{eq:22}) and get an almost periodic word
\begin{equation}
	\label{eq:30}
	W_{\a} = \blv_{1} \blv_{2}\underbrace{ \blv_{1/2}\blv_{1/2}}_{s}\blv_3 \blv_4 \dots \blv_{k-1} \blv_{k} \underbrace{\blv_{1/2}\blv_{1/2}}_{r} \blv_3 \blv_4 \dots,
\end{equation}
where $s \geq 0, r \geq 2$ are  arbitrary integers. Here $\blv_{1/2}$ stands either for $\blv_1$ or $\blv_2$, because $\pr^{-1}_{\bbeta_1}(\blu_1) = \{\blv_1, \blv_2\}$.
\import{./}{pic_8.tex}

We first show that we cannot have
\begin{equation*}
	W_{\a} = \blv_1 \blv_2 \blv_1 \dots.
\end{equation*}
Note that by (\ref{eq:30}), $t_1(\bbeta_1) = t_{s+2}(\a)$, because it is a moment of change of value of $\blv_{\a}(t)$ from $\blv_{1/2}$ to $\blv_3$, and, therefore, is a moment of change of value of $\blv_{\bbeta_1}(t)$ from $\blu_{1}$ to $\blu_2$ (see Figure \ref{fig:per}).

By \ref{it1}, there are again $m = \frac{k(k-1)}{2}$ functions $\psi_{j_1}, \dots, \psi_{j_m}$ which are continuous at $t_{s+2}(\a)$. We denote the corresponding subcollection of numbers to these functions by $\bbeta_2 = (\alpha_{j_1}, \dots, \alpha_{j_m})$. This subcollection is different from $\bbeta_1$, because the value of $\blv_{\bbeta_1}(t)$ changes at $t_{s+2}(\a)$. By the same argument as for $\bbeta_1$, we can apply the induction hypothesis to $\bbeta_2$. So, there exists $T(\bbeta_2)$ such that for all $t'_0 \geq T(\bbeta_2)$ we get a $(k-1)$-periodic sequence $\{\blv_{\bbeta_2}(t_i(\bbeta_2, t'_0))\}_{i=0}^\infty$, where values $\blv_{\bbeta_2}(t_{i}(\bbeta_2, t'_0)),$ $0 \leq i \leq k-2$ are all different. Again, since we are proving the existence of $T = T(\a)$, if $T(\bbeta_2) \geq T$, we can set $T = T(\bbeta_2)$, chose a new arbitrary $t_0 \geq T$, and repeat the arguments above for this new values.

%Here, the sequence $\{t_i(\bbeta_2)\}_{i=0}^\infty$ is different from $\{t_i(\bbeta_1)\}_{i=0}^\infty$.
We denote
\begin{equation*}
	\blw_{i+1} = \blv_{\bbeta_2}(t_i(\bbeta_2, t_0)), \ 0 \leq i \leq k-2,
\end{equation*}
so, in the word notation, we have a $(k-1)$-periodic word
\begin{equation}
	\label{eq:25}
	W_{\bbeta_2} = \overline{ \blw_1 \blw_2  \dots \blw_{k-2} \blw_{k-1}\vphantom{T}}.
\end{equation} 
The functions $\psi_{j_1}, \dots, \psi_{j_m}$ are continuous at $t_{s+2}(\a)$, so, by the same argument as for $\bbeta_1$ we have
\begin{equation*}
	\blv_{\bbeta_2}(t_{s+1}(\a)) = \blv_{\bbeta_2}(t_{s+2}(\a)).
\end{equation*}
In the interval $[t_0, t_{s+1}]$ only two values $\blv_1, \blv_2$ of $\blv_{\a}(t)$ appear. So, in the same interval, no more than two values of $\blv_{\bbeta_2}(t)$ can appear, that is,
$\blv_{\bbeta}(t_{s+1}(\a))$ is either $\blw_1$ or $\blw_2$. We denote this as
\begin{equation}
	\label{eq:14}
	\blw_{1/2} = \blv_{\bbeta}(t_{s+1}(\a)) = \blv_{\bbeta}(t_{s+2}(\a)).
\end{equation}
Since $t_{s+2}(\a)$ is a moment of change of value of $\blv_{\a}(t)$ from $\blv_{1/2}$ to $\blv_3$,  $$\blv_{1/2}=\blv_{\a}(t_{s+1}(\a)) \neq \blv_{\a}(t_{s+2}(\a)) = \blv_3.$$ By the same argument as for (\ref{eq:23}), we get 
\begin{equation}
	\label{eq:34}
	\{\blv_{1/2}, \blv_3\} \subset \pr^{-1}_{\bbeta_2}(\blw_{1/2}).
\end{equation}
In the same way as for $\bbeta_1$, analogously to (\ref{eq:31}) and (\ref{eq:32}), we have
\begin{align}
	\label{eq:42}
	&|\pr^{-1}_{\bbeta_2}(\blw_{1/2})| = 2, \\
	\label{eq:43}
	&|\pr^{-1}_{\bbeta_2}(\blw_i)| = 1, \ 1 \leq i \leq k-1, \ \blw_i \neq \blw_{1/2}, \\
	\label{eq:44}
	&\pr^{-1}_{\bbeta_2}(\blw_i) \cap \pr^{-1}_{\bbeta_2}(\blw_j) = \emptyset, \ i \neq j.
\end{align}
By (\ref{eq:25}) and because $k \geq 4$, we cannot have 
\begin{equation*}
	W = \blw_1 \blw_2 \blw_1 \blw_3 \dots.
\end{equation*}
So, if in (\ref{eq:14}) we have $\blw_{1/2} = \blw_1$, it follows that
\begin{equation*}
	\blv_{\bbeta_2}(t) = \blw_{1} \text{ for all } t \in \big[t_0, t_{s+2}(\a)\big].
\end{equation*}
Therefore,
\begin{equation*}
	\{\blv_{1}, \blv_2, \blv_3\} \subset \pr^{-1}_{\bbeta_2}(\blw_1),
\end{equation*}
which contradicts (\ref{eq:42}).
Thus, $\blw_{1/2} = \blw_2$.
So, we can rewrite (\ref{eq:42}), (\ref{eq:43}) as
\begin{align}
	\label{eq:45}
	&|\pr^{-1}_{\bbeta_2}(\blw_{2})| = 2, \\
	\label{eq:46}
	&|\pr^{-1}_{\bbeta_2}(\blw_i)| = 1, \ 1 \leq i \leq k-1, \ i \neq 2.
\end{align}
And since
\begin{equation*}
	\{\blv_1\} \subset \pr^{-1}_{\bbeta_2}(\blw_1), 
\end{equation*}
by (\ref{eq:34}), (\ref{eq:44}), (\ref{eq:45}) and (\ref{eq:46}), we have
\begin{align}
	\label{eq:36}
	&\pr^{-1}_{\bbeta_2}(\blw_1) = \blv_1, \\
	\label{eq:37}
	&\pr^{-1}_{\bbeta_2}(\blw_2) = \{\blv_2, \blv_3\}, \\ \label{eq:38}
	&\pr^{-1}_{\bbeta_2}(\blw_i) = \blv_{i+1}, \ 3 \leq i \leq k-1.
\end{align}
Assume that we have
\begin{equation*}
	W_{\a} = \blv_1 \blv_2 \blv_1 \dots.
\end{equation*}
For the subcollection $\bbeta_2$ by (\ref{eq:36})--(\ref{eq:38}) it corresponds to
\begin{equation*}
	W_{\bbeta_2} = \blw_1 \blw_2 \blw_1 \dots,
\end{equation*}
which is a contradiction with (\ref{eq:25}), because  $k \geq 4$ and we have
\begin{equation*}
	W_{\bbeta_2} = \blw_1 \blw_2 \blw_3 \dots.
\end{equation*}
The same argument shows that the entire word $W_{\a}$ does not contain subword $\blv_1\blv_2\blv_1$, not just the beginning.

We are left to show that in (\ref{eq:30}) there are no subwords $\blv_{k} \blv_2 \blv_1$. If there is such a subword in $W_{\a}$, by (\ref{eq:36})--(\ref{eq:38}) it corresponds to
$\blw_{k-1} \blw_2 \blw_1 ,$ which is a contradiction with (\ref{eq:25}).
This completes the proof of \ref{it2}.

\subsection{Proof of \ref{it3} - \ref{it6}} We now prove \ref{it3} - \ref{it6} simultaneously by induction, starting with $k=2$.

\textbf{The base case.}
Let $k=2$ and $n=3$. Let $\a = (\alpha_1, \alpha_2, \alpha_3)$ be three pairwise independent numbers with $\knm(\a) = 2$. Let $\blv_1, \, \blv_2$ be two distinct values of $\blv_{\a}(t)$ with $\blv_{\a}(t_0) = \blv_1, \ \blv_{\a}(t_1) = \blv_2$.
\import{./}{pic_12.tex}
Recall that we enumerate functions as $\psi_1, \psi_{1,2}, \psi_2$, so that at $t_0$ we have
\begin{equation*}
	\psi_1(t_0) > \psi_{1,2}(t_0) > \psi_{2}(t_0).
\end{equation*}
This case is illustrated in Figure \ref{fig:ex3}.
By Kan-Moshchevitin Theorem \ref{th:KM}, any two functions have to interchange, so two vectors $\blv_1,  \blv_2$ are reversed with respect to each other, that is, if $\blv_1 = (v_1, v_2, v_3)$, then $\blv_2 = (v_3, v_2, v_1)$. So, 
\begin{equation*}
	\psi_2(t_1) > \psi_{1,2}(t_1) > \psi_1(t_1).
\end{equation*}
Thus, $\psi_1$ and $\psi_{1,2}$ have a jump at $t_1$ (see Figure \ref{fig:ex3}). Since $\blv_{\a}(t_{2i}) = \blv_1, \ i \geq 0$, we have
\begin{equation*}
	\psi_1(t_{2i}) > \psi_{1,2}(t_{2i}) > \psi_{2}(t_{2i}), \ i \geq 0.
\end{equation*}
And because $\blv_{\a}(t_{2i+1}) = \blv_2, \ i \geq 0$, we have
\begin{equation*}
	\psi_1(t_{2i+1}) > \psi_{1,2}(t_{2i+1}) > \psi_{2}(t_{2i+1}), \ i \geq 0.
\end{equation*}
So, functions $\psi_1$ and $\psi_{1,2}$ have  jumps at $t_{2i+1}, i \geq 0$, and functions $\psi_2$ and $\psi_{1,2}$ have  jumps at $t_{2i}, i \geq 1$. 

Functions $\psi_{2,2}$ and $\psi_{1,2}$ have a simultaneous jump at $t_2$. Let $t'$ denote the previous jump of $\psi_{1,2}$. By Lemma \ref{lm:MM}, since $$\psi_{2,2}(t_2-1) > \psi_{1,2}(t_2-1),$$ we have $$\psi_{1,2}(t'-1) > \psi_{2,2}(t'-1).$$ So, $t' = t_1$, thus
function $\psi_{1,2}$ is continuous between $t_1$ and $t_2$. The same argument works for all $i \geq 0$, hence the function $\psi_{1,2}$ is continuous between $t_i$ and $t_{i+1}$ for all $i \geq 0$. This completes the proof for $k=2, \, n = 3$.

\textbf{Step of induction}. Let $k \geq 3$. As the induction hypothesis, we assume that for any $\frac{k(k-1)}{2}$-tuple $\bbeta$ of pairwise independent numbers with $\knm(\bbeta) = k-1$ items \ref{it3} -- \ref{it6} hold. We will prove that for a $\frac{k(k+1)}{2}$-tuple  $\a$ with $\knm(\a) = k$ items \ref{it3} -- \ref{it6} also hold. The case $k=3, n=6$ is illustrated in Figure \ref{fig2}.

By \ref{it2}, $\blv_{\a}(t_{i}(\a)), 0 \leq i \leq k-1$ are $k$ distinct values of $\blv_{\a}(t)$ that change periodically. We denote
\begin{equation*}
	\blv_i = \blv_{\a}(t_{i-1}(\a)),  \, 1 \leq i \leq k.
\end{equation*}
By periodicity, for all $s \geq 0 $ we have
\begin{equation*}
	\blv_{\a}(t_{i-1 + sk}(\a)) = \blv_i,  \, 1 \leq i \leq k.
\end{equation*}
By \ref{it1} applied to $\a$ there are $m = \frac{k(k-1)}{2}$ functions $\psi_{i_1}, \dots, \psi_{i_m}$ which are continuous at $t_{k}(\a)$. We denote the  subcollection of numbers corresponding to  functions $\psi_{i_1}, \dots, \psi_{i_m}$ by $\bbeta_1 = (\alpha_{i_1}, \dots, \alpha_{i_m})$. For example, in case $k=3$ (see Figure \ref{fig2}), we have $t_k(\a) = t_3(\a)$, and functions that are continuous at $t_3(\a)$ are functions $\psi_1, \, \psi_3, \, \psi_4$.

As in the proof of \ref{it2}, the vector of these functions does not change at $t_k(\a)$, so 
\begin{equation}
	\label{eq:47}
	\blv_{\bbeta_1}(t_{k-1}(\a)) =\blv_{\bbeta_1}(t_{k}(\a)),
\end{equation}
and, by the same arguments as before,
\begin{equation*}
	\knm(\bbeta_1) = k-1.
\end{equation*}
By \ref{it2}, $\blv_{\bbeta_1}(t_{i}(\bbeta_1)), \, 0 \leq i \leq k-1$ are $k-1$ distinct values of $\blv_{\bbeta_1}(t)$ that change periodically. We denote
\begin{equation}
	\label{eq:48}
	\blu_{i} = \blv_{\bbeta_1}(t_{i-1}(\bbeta_1)), \ 1 \leq i \leq k-1.
\end{equation}
By periodicity, we have
\begin{equation}
	\label{eq:54}
	\blv_{\bbeta_1}(t_{k-1}(\bbeta_1)) = \blu_1.
\end{equation}
Because the point of change of value of $\blv_{\bbeta_1}(t)$ is also the point of change of value of $\blv_{\a}(t)$, it follows that
\begin{equation}
	\label{eq:55}
	t_{i}(\bbeta_1) = t_{i}(\a), \ 0 \leq i \leq k-1.
\end{equation}
However, functions $\bbeta_1$ are continuous at $t_{k}(\a)$, so
\begin{equation*}
	t_k(\bbeta_1) = t_{k+1}(\a).
\end{equation*}
\import{./}{pic_2.tex}
By (\ref{eq:47}) -- (\ref{eq:55}), functions $\bbeta_1$ have the same vector at $t_0, t_{k-1}(\a)$ and $t_k(\a)$, namely,
\begin{equation}
	\label{eq:57}
	\blu_1 = \blv_{\bbeta_1}(t_0) =\blv_{\bbeta_1}(t_{k-1}(\a)) = \blv_{\bbeta_1}(t_k(\a))
\end{equation}
(here in the second equality we use $t_{k-1}(\bbeta_1) = t_{k-1}(\a)$ from (\ref{eq:55})).

Now it would be convenient to use our special enumeration (\ref{int:eq10}).
We enumerate $m=\frac{k(k-1)}{2}$ functions  $\psi_{i_1}, \dots, \psi_{i_m}$ with indices $(i,j), \, 1 \leq i \leq k-1, \ i \leq j \leq k-1$, in such a way that at $t_0$ they are ordered as
\begin{multline}
	\label{eq:49}
	\psi_{1,1}(t_0) > \psi_{1,k-1}(t_0) > \psi_{1,k-2}(t_0) > \dots > \psi_{1,2}(t_0) > \\ >  \psi_{2,2}(t_0) > \psi_{2,k-1}(t_0) > \psi_{2,k-2}(t_0) > \dots > \psi_{2,3}(t_0) > \\ > \psi_{3,3}(t_0)  > \dots > \psi_{k-2, k-2}(t_0) > \psi_{k-2, k-1}(t_0) > \psi_{k-1,k-1}(t_0),
\end{multline}
In Figure \ref{fig2} these functions are the functions $\psi_{1,1}, \, \psi_{1,2}, \, \psi_{2,2}$.

By induction hypothesis items \ref{it3} -- \ref{it6} hold for functions $\bbeta_1$. So,
\begin{enumerate}[label=(\alph*)]
	\item \label{it3:pr} Each function $\psi_{i,i}, \ 1 \leq i \leq k-1$ has jumps at $\{t_{i+sk}(\bbeta_1)\}_{s=0}^\infty$ and is continuous at all other $\{t_j(\bbeta_1)\}_{j=1}^\infty$.
	\item \label{it4:pr} Each function $\psi_{i,j}, \ 1 \leq i \leq k, \ i < j \leq k$ has jumps at $\{t_{i+sk}(\bbeta_1)\}_{s=0}^\infty$ and $\{t_{j+sk}(\bbeta_1)\}_{s=0}^\infty$, and is continuous everywhere between these jumps.
\end{enumerate}

Let $\psi_{j_1}, \dots, \psi_{j_k}$ be the remaining $k$ functions that have a jump at $t_k(\a)$ (in Figure \ref{fig2} these functions are $\psi_3, \, \psi_5, \, \psi_6$). By \ref{it1}, at each $t_{i}(\bbeta_1)$ exactly $k-1$ functions from $\bbeta_1$ have a jump. So, at each $t_{i}(\bbeta_1) = t_i(\a), 1 \leq i \leq k-1$ only one function from $\psi_{j_1}, \dots, \psi_{j_k}$ has a jump. 

Let $\psi_{j_r}$ be the largest function at $t_{k-1}(\a)$ from the subcollection $\psi_{j_1}, \dots, \psi_{j_k}$ (function $\psi_6$ in Figure \ref{fig2}). By Kan-Moshchevitin Theorem \ref{th:KM}, each of the remaining $k-1$ functions is larger than $\psi_{j_r}$ at some $t \in \big[t_0, t_{k-1}(\a)\big)$. So, for each function $\psi_{j_l} \neq \psi_{j_r}$ there exist $i, \ 1 \leq i \leq k-1$, such that $\psi_{j_l}$ has a jump at $t_{i}(\a)$ and
\begin{equation}
	\label{eq:59}
	\psi_{j_l}\big(t_i(\a) - 1\big) > \psi_{j_r}\big(t_i(\a) - 1\big)
\end{equation}
(recall that by $\psi_{i}(t-1)$ we denote the value of $\psi_{i}$ just before $t$).

Since at each $t_{i}(\a), \, 1 \leq i \leq k-1$ exactly one function from $\psi_{j_1}, \dots, \psi_{j_k}$ has a jump, it follows that for each $\psi_{j_s} \neq \psi_{j_r}$ there exists only one $t_{i}(\a), \, 1 \leq i \leq k,$ such that $\psi_{j_s}$ has a jump at $t_i(\a)$, and for different $\psi_{j_s}$, $\psi_{j_l}$ such $t_s(\a), t_l(\a)$ are different. It also follows that $\psi_{j_r}$ is continuous at all $t_{i}(\a), \, 1 \leq i \leq k-1$.

We enumerate functions $\psi_{j_1}, \dots, \psi_{j_k}$ with indices $(i,k),  \ 1 \leq i \leq k$ in such a way that 
$\psi_{i,k}$ has jumps at $t_i(\a), 1 \leq i \leq k-1,$ and $t_k(\a)$, with $\psi_{k,k}$ denoting $\psi_{j_r}$ which has a jump only at $t_{k}(\a)$. By construction $\psi_{k,k} = \psi_{j_r}$ is the largest function at $t_{k-1}(\a)$, so
\begin{equation*}
	\psi_{k,k}(t_{k-1}(\a)) > \psi_{i,k}(t_{k-1}(\a)), \ 1 \leq i \leq k-1.
\end{equation*}
By (\ref{eq:59}) we have
\begin{equation}
	\label{eq:60}
	\psi_{i,k}\big(t_i(\a) - 1\big) > \psi_{k,k}\big(t_i(\a) - 1\big).
\end{equation}

In Figure \ref{fig2} function $\psi_{3,3} = \psi_6$ is the largest function at $t_2(\a)$. Functions $\psi_{1,3} = \psi_{2}$ and $\psi_{2,3} = \psi_{5}$ have jumps at $t_1(\a)$ and $t_2(\a)$ respectively.

By Lemma \ref{lm:MM}, functions $\psi_{i,k}, \ 1 \leq i \leq k-1$ are continuous between jumps $t_i(\a)$ and $t_k(\a)$ (we use the same argument as at the conclusion of the proof of the base case).

We have considered $t \in [t_0, t_{k}].$ Now we explore periodicity.

The same $k$ functions have a simultaneous jump at each $t_{ks}(\a), s \in \N_0$. By the same argument and periodicity of the sequence $\{\blv_{\a}(t_{i}(\a))\}_{i=0}^{\infty}$ we get that functions $\psi_{i,k}, \ 1 \leq i \leq k-1$ have jumps at $t_{i+sk}(\a)$ and $t_{ks}(\a)$, and are continuous between them. The function $\psi_{k,k}$ has jumps at $t_{ks}(\a)$, and is continuous at all other $t_{i}(\a)$. Together with \ref{it3:pr}, \ref{it4:pr} this proves \ref{it3}, \ref{it4}.

Now we prove \ref{it5}.
For all $1 \leq i \leq k-2$ we apply Lemma \ref{lm:3_func} to functions $\psi_{k-i, k}, \, \psi_{k-i-1,k}, \, \psi_{k-i-1, k-i}$, with $\psi_{\alpha} = \psi_{k-i-1,k}, \ \psi_{\beta} = \psi_{k-i,k}$ and $\psi_{\gamma} = \psi_{k-i-1,k-i}$. Functions $\psi_{k-i, k}, \, \psi_{k-i-1,k}$ have simultaneous jumps at $t_{k+kj}(\a), \, j \in \N_0$. Functions $\psi_{k-i, k}, \, \psi_{k-i-1,k-i}$ have simultaneous jumps at $t_{k-i+kj}(\a), \, j \in \N_0$. Functions $\psi_{k-i-1, k}, \, \psi_{k-i-1,k-i}$ have simultaneous jumps at $t_{k-i-1+kj}(\a), \, j \in \N_0$. For all $1 \leq i \leq k-2$ we get
\begin{equation*}
	\psi_{k-i, k}(t)> \psi_{k-i-1, k}(t), \text{ for all } t \in \big[t_{k-i}(\a), t_{k}(\a)\big),
\end{equation*}
In particular, we have
\begin{equation*}
	\psi_{k-i, k}\big(t_{k-1}(\a)\big) > \psi_{k-i-1, k}\big(t_{k-1}(\a)\big), \text{ for all } 1 \leq i \leq k-2.
\end{equation*}
Hence,
\begin{equation}
	\label{eq:50}
	\psi_{k,k}\big(t_{k-1}(\a)\big) > \psi_{k-1, k}\big(t_{k-1}(\a)\big) > \psi_{k-2, k}\big(t_{k-1}(\a)\big) > \dots > \psi_{1, k}\big(t_{k-1}(\a)\big).
\end{equation}
For example, in Figure \ref{fig2} we have
\begin{equation*}
	\psi_{3,3}\big(t_{2}(\a)\big) > \psi_{2, 3}\big(t_{2}(\a)\big) > \psi_{1, 3}\big(t_{2}(\a)\big).
\end{equation*}
Functions $\psi_{1,k}$ and $\psi_{1,1}$ are both continuous at $t_{k-1}(\a)$, and thus, by \ref{it1}, belong to the subcollection of $\frac{k(k-1)}{2}$ functions that are continuous at $t_{k-1}(\a)$. We denote the corresponding subcollection of numbers by $\bbeta_2$. 

In Figure \ref{fig2} the functions that are continuous at $t_2(\a)$ are the functions $\psi_{1,1}, \psi_{1,3}$ and $\psi_{3,3}$.

We apply the induction hypothesis for \ref{it5} to $\bbeta_2$. The functions that have a jump at each $t_i(\bbeta_2)$ are the $k-1$ largest functions at $t_{i-1}(\bbeta_2)$. So, since $\psi_{1,k}$ has a jump at $t_k(\a)$ and $\psi_{1,1}$ is continuous at $t_k(\a)$, we get
\begin{equation}
	\label{eq:56}
	\psi_{1,k}(t_{k-1}(\a)) > \psi_{1,1}(t_{k-1}(\a)) .
\end{equation}
By (\ref{eq:57}), inequality (\ref{eq:49}) also holds for $t_{k-1}(\a)$, so
\begin{multline}
	\label{eq:58}
	\psi_{1,1}\big(t_{k-1}(\a)\big) > \psi_{1,k-1}\big(t_{k-1}(\a)\big) > \psi_{1,k-2}\big(t_{k-1}(\a)\big) > \dots > \psi_{1,2}\big(t_{k-1}(\a)\big) > \\ >  \psi_{2,2}\big(t_{k-1}(\a)\big) > \psi_{2,k-1}\big(t_{k-1}(\a)\big) > \psi_{2,k-2}\big(t_{k-1}(\a)\big) > \dots > \psi_{2,3}\big(t_{k-1}(\a)\big) > \\ > \psi_{3,3}\big(t_{k-1}(\a)\big)  > \dots > \psi_{k-2, k-2}\big(t_{k-1}(\a)\big) > \psi_{k-2, k-1}\big(t_{k-1}(\a)\big) > \psi_{k-1,k-1}\big(t_{k-1}(\a)\big),
\end{multline}
By (\ref{eq:50}), (\ref{eq:56}) and (\ref{eq:58}), we get
\begin{multline}
	\label{eq:52}
	\psi_{k,k}\big(t_{k-1}(\a)\big) > \psi_{k-1, k}\big(t_{k-1}(\a)\big) > \psi_{k-2, k}\big(t_{k-1}(\a)\big) > \dots > \psi_{1, k}\big(t_{k-1}(\a)\big) > \\ >
	\psi_{1,1}\big(t_{k-1}(\a)\big) > \psi_{1,k-1}\big(t_{k-1}(\a)\big) > \psi_{1,k-2}\big(t_{k-1}(\a)\big) > \dots > \psi_{1,2}\big(t_{k-1}(\a)\big) > \\  > \psi_{2,2}\big(t_{k-1}(\a)\big)  > \dots > \psi_{k-2, k-2}\big(t_{k-1}(\a)\big)> \psi_{k-2, k-1}\big(t_{k-1}(\a)\big) > \psi_{k-1,k-1}\big(t_{k-1}(\a)\big).
\end{multline}
By construction, functions $\psi_{k,k}, \psi_{1,k}, \dots, \psi_{1,2}$ have a jump at $t_{k}(\a)$. So, \ref{eq:52} shows that the largest $k$ functions at $t_{k-1}(\a)$ have a jump at $t_k(\a)$, and $\psi_{k,k}$ is the largest function. 

The same argument works for all $t_{sk}(\a), \, s \in \N_0$. Since $t_0$ is arbitrary, we can set $t_0 = t_i(\a)$. Then the new $t_k(\a)$ is $t_{k+i}(\a)$, and by the same steps we get that the functions that have a jump at $t_{k+i}(\a)$ are the $k$ largest functions at $t_{k+i - 1}(\a)$, and $\psi_{i,i}$ is the largest function at $t_{k+i - 1}(\a)$. By periodicity, it follows that the largest $k$ functions at $t_{i-1}(\a)$ have a jump at $t_i(\a)$, with $\psi_{i,i}$ being the largest function $t_{i-1}(\a)$ for all $i \geq 1$. This \mbox{proves \ref{it5}.}

Finally, we prove \ref{it6}.
For all $ 1 \leq i \leq k-2$, we apply Lemma \ref{lm:3_func} to $\psi_{\alpha} = \psi_{i, k-1}, \, \psi_{\beta} = \psi_{i, k}$ and $\psi_{\gamma} = \psi_{k-1, k}$. We obtain
\begin{equation*}
	\psi_{i, k}(t) > \psi_{i, k-1}(t), \text{ for all } t \in \big[t_k(\a), t_{i+k}(\a)\big).
\end{equation*}
In particular,
\begin{equation}
	\label{eq:51}
	\psi_{i, k}\big(t_{k}(\a)\big) > \psi_{i, k-1}\big(t_{k}(\a)\big), \text{ for all } 1 \leq i \leq k-2.
\end{equation}
For all $1 \leq i \leq k-1$, by \ref{it5}, function $\psi_{i,i}$ is the largest function at $t_{i-1+k}(\a)$, so
\begin{equation*}
	\psi_{i,i}\big(t_{i-1+k}(\a)\big) > \psi_{i,k}\big(t_{i-1+k}(\a)\big).
\end{equation*}
Because functions $\psi_{i, k}$ and $\psi_{i,i}$ are continuous at all $t_{k+j}, \, 1 \leq j \leq i-1$, at $t_{k}(\a)$ we also have
\begin{equation*}
	\psi_{i,i}\big(t_{k}(\a)\big) > \psi_{i,k}\big(t_{k}(\a)\big) \text{ for all } 1 \leq i \leq k-1.
\end{equation*} 
By \ref{eq:60} we have
\begin{equation*}
	\psi_{k-1,k}(t_{k-1}(\a)-1) > \psi_{k, k}(t_{k-1}(\a)-1).
\end{equation*}
Because functions $\psi_{k-1,k}$ and $\psi_{k, k}$ are continuous at all $t_{j}(\a), \, 1 \leq j \leq k-2$, we have
\begin{equation}
	\label{eq:61}
	\psi_{k-1,k}(t_0) > \psi_{k, k}(t_0).
\end{equation}
Since $\blv_{\a}(t_0) = \blv_{\a}\big(t_{k}(\a)\big)$, we have
\begin{equation}
	\label{eq:62}
	\psi_{k-1,k}\big(t_{k}(\a)\big) > \psi_{k, k}\big(t_{k}(\a)\big).
\end{equation}
So, from (\ref{eq:51}), (\ref{eq:61}) and (\ref{eq:62}), we get
\begin{multline}
	\label{eq:63}
	\psi_{1,1}(t_{k}(\a)) > \psi_{1,k}(t_{k}(\a)) > \psi_{1,k-1}(t_{k}(\a)) > \dots > \psi_{1,2}(t_{k}(\a)) > \\ >  \psi_{2,2}(t_{k}(\a)) > \psi_{2,k}(t_{k}(\a)) > \psi_{2,k-1}(t_{k}(\a)) > \dots > \psi_{2,3}(t_{k}(\a)) > \\ > \psi_{3,3}(t_{k}(\a))  > \dots > \psi_{k-1}(t_{k}(\a)) > \psi_{k-1, k}(t_{k}(\a)) > \psi_{k,k}(t_{k}(\a)).
\end{multline}
From (\ref{eq:52}) and (\ref{eq:63}), by Remark \ref{rm:pi}, we have
\begin{equation*}
	\blv_{\alpha}(t_{k}(\a)) = \pi\bigl( \blv_{\alpha}(t_{k-1}(\a))\bigr).
\end{equation*}
Since $t_0$ was arbitrary, we can set $t_0 = t_i(\a), \ i \geq 0$. Then, the new $t_k(\a)$ is $t_{k+i}(\a)$, and by repeating the same steps, we get
\begin{equation}
	\label{eq:53}
	\blv_{\alpha}(t_{k+i}(\a)) = \pi\bigl( \blv_{\alpha}(t_{k+i-1}(\a))\bigr), \ i \geq 0.
\end{equation}
Since $\blv_{\alpha}(t_{k+i}(\a)) = \blv_{\alpha}(t_{i}(\a))$ for all $0 \leq i \leq k-1$, equality (\ref{eq:52}) also holds for $t_{i}(\a), \ 0 \leq i \leq k-1.$
Thus, we have
\begin{equation*}
	\blv_{\alpha}(t_{i}(\a)) = \pi\bigl( \blv_{\alpha}(t_{i-1}(\a))\bigr), \ i \geq 1.
\end{equation*}
This completes the proof of Theorem \ref{th:main}.
\begin{comment}

\begin{remark}
	\label{rm:3}
	In the proof of \ref{it2} we check the base case for $k = 3$. It is necessary, because the proof of the inductive step does not work when $k = 3$. Here we provide an example. 
\end{remark}

\textit{Left to do:}
Acknowledgment, example 3.7, remark 5.1, figure 1?, figure 7 add a function, remark 3.4.
\end{comment}
\subsection*{Acknowledgements}
	The author is grateful to Nikolay Moshchevitin for formulating this problem and for many helpful discussions. The author also thanks Erez Nesharim for reading this paper and pointing out several mistakes.

%% file: pic_10.tex
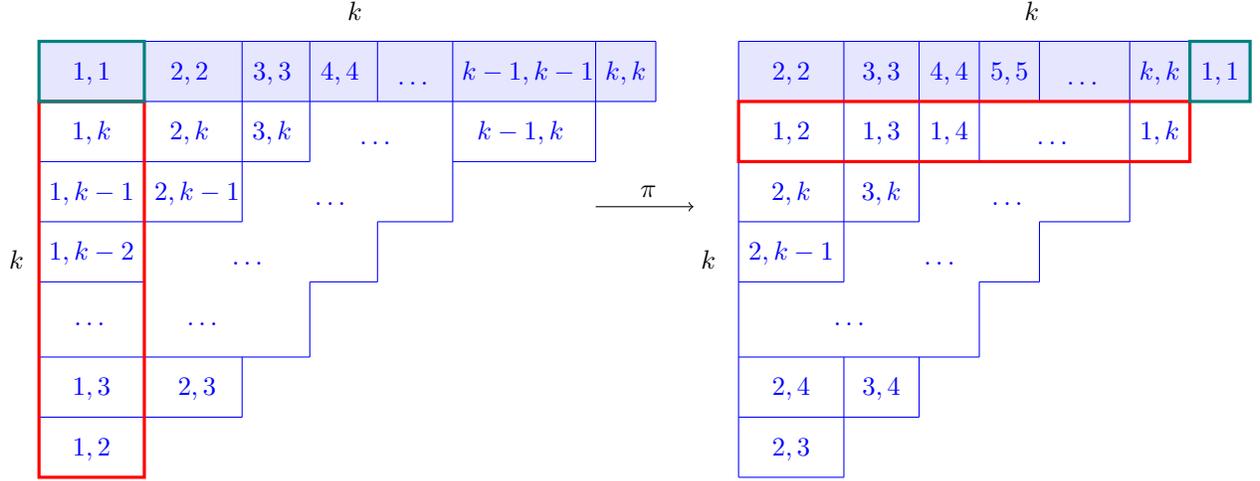
\begin{figure}[t!]
	\centering
	\begin{tikzpicture}	
		\fill[color=blue!10] (-0.7,3.2) rectangle (7.5,4);
		\draw[blue!100] (-0.7, -1.8)  -- (-0.7, 4);
		\draw[blue!100] (-0.7, 4)  -- (7.5, 4);
		\draw[blue!100] (-0.7, 3.2) node[blue, right=0.7cm, above=0.1cm] {$1,1$}  node[blue, right=2cm, above=0.1cm] {$2,2$} node[blue, right=3.1cm, above=0.1cm] {$3,3$} node[blue, right=4cm, above=0.1cm] {$4,4$} node[blue, right=5cm, above=0.1cm] {$\dots$} node[blue, right=6.5cm, above=0.1cm] {$k-1,k-1$} node[blue, right=7.8cm, above=0.1cm] {$k,k$} -- (7.5, 3.2);
		\draw[blue!100] (0.7, -1.8)  -- (0.7, 4);
		\draw[blue!100] (2, 1.6)  -- (2, 4);
		\draw[blue!100] (2.9, 2.4)  -- (2.9, 4);
		\draw[blue!100] (3.8, 3.2)  -- (3.8, 4);
		\draw[blue!100] (4.8, 3.2)  -- (4.8, 4);
		\draw[blue!100] (6.7, 3.2)  -- (6.7, 4);
		\draw[blue!100] (7.5, 3.2)  -- (7.5, 4);
		
		\draw[blue!100] (2, -1)  -- (2, -0.2);
		\draw[blue!100] (2.9, -0.2)  -- (2.9, 0.8);
		\draw[blue!100] (2.9, 0.8)  -- (3.8, 0.8);
		\draw[blue!100] (3.8, 0.8)  -- (3.8, 1.6);
		\draw[blue!100] (3.8, 1.6)  -- (4.8, 1.6);
		\draw[blue!100] (4.8, 1.6)  -- (4.8, 3.2);
		\draw[blue!100] (4.8, 2.4) node[blue, right=0.9cm, above=0.1cm] {$k-1,k$} -- (6.7, 2.4);
		\draw[blue!100] (6.7, 2.4)  -- (6.7, 3.2);
		
		\draw[blue!100] (-0.7, 2.4) node[blue, right=0.7cm, above=0.1cm] {$1,k$}  node[blue, right=2cm, above=0.1cm] {$2,k$} node[blue, right=3.1cm, above=0.1cm] {$3,k$} node[blue, right=4.5cm, above=0.1cm] {$\dots$}  -- (2.9, 2.4);
		\draw[blue!100] (-0.7, 1.6) node[blue, right=0.7cm, above=0.1cm] {$1,k-1$}  node[blue, right=2.1cm, above=0.1cm] {$2,k-1$} node[blue, right=3.9cm, above=0.1cm] {$\dots$}  -- (2, 1.6);
		\draw[blue!100] (-0.7, 0.8) node[blue, right=0.7cm, above=0.1cm] {$1,k-2$}  node[blue, right=2.8cm, above=0.1cm] {$\dots$} -- (0.7, 0.8);
		\draw[blue!100] (-0.7, -0.2) node[blue, right=0.7cm, above=0.3cm] {$\dots$} node[blue, right=2.2cm, above=0.3cm] {$\dots$}   -- (2.9, -0.2);
		\draw[blue!100] (-0.7, -1) node[blue, right=0.7cm, above=0.1cm] {$1,3$} node[blue, right=2.1cm, above=0.1cm] {$2,3$}  -- (2, -1);
		\draw[blue!100] (-0.7, -1.8) node[blue, right=0.7cm, above=0.1cm] {$1,2$}   -- (0.7, -1.8);
		
		\fill[color=blue!10] (8.6,3.2) rectangle (15.4,4);
		\draw[blue!100] (8.6, -1.8)  -- (8.6, 4);
		\draw[blue!100] (8.6, 4)  -- (15.4, 4);
		\draw[blue!100] (8.6, 3.2) node[blue, right=0.7cm, above=0.1cm] {$2,2$}  
		node[blue, right=1.9cm, above=0.1cm] {$3,3$} 
		node[blue, right=2.8cm, above=0.1cm] {$4,4$} 
		node[blue, right=3.6cm, above=0.1cm] {$5,5$} 
		node[blue, right=4.6cm, above=0.1cm] {$\dots$} 
		node[blue, right=5.6cm, above=0.1cm] {$k,k$} 
		node[blue, right=6.4cm, above=0.1cm] {$1,1$} -- (15.4, 3.2);
		\draw[blue!100] (10, -1.8)  -- (10, -0.2);
		\draw[blue!100] (10, 0.8)  -- (10, 4);
		\draw[blue!100] (11, 1.6)  -- (11, 4);
		\draw[blue!100] (11.8, 2.4)  -- (11.8, 4);
		\draw[blue!100] (12.6, 3.2)  -- (12.6, 4);
		\draw[blue!100] (13.8, 3.2)  -- (13.8, 4);
		\draw[blue!100] (14.6, 3.2)  -- (14.6, 4);
		\draw[blue!100] (15.4, 3.2)  -- (15.4, 4);
		
		\draw[blue!100] (11, -1)  -- (11, -0.2);
		\draw[blue!100] (11.8, -0.2)  -- (11.8, 0.8);
		\draw[blue!100] (11.8, 0.8)  -- (12.6, 0.8);
		\draw[blue!100] (12.6, 0.8)  -- (12.6, 1.6);
		\draw[blue!100] (12.6, 1.6)  -- (13.8, 1.6);
		\draw[blue!100] (13.8, 1.6)  -- (13.8, 3.2);
		\draw[blue!100] (13.8, 2.4) node[blue, right=0.4cm, above=0.1cm] {$1,k$} -- (14.6, 2.4);
		\draw[blue!100] (14.6, 2.4)  -- (14.6, 3.2);
		
		\draw[blue!100] (8.6, 2.4) node[blue, right=0.7cm, above=0.1cm] {$1,2$}  
		node[blue, right=1.9cm, above=0.1cm] {$1,3$} 
		node[blue, right=2.8cm, above=0.1cm] {$1,4$} 
		node[blue, right=4.2cm, above=0.1cm] {$\dots$}  
		-- (11.8, 2.4);
		\draw[blue!100] (8.6, 1.6) node[blue, right=0.7cm, above=0.1cm] {$2,k$}  
		node[blue, right=1.9cm, above=0.1cm] {$3,k$} 
		node[blue, right=3.6cm, above=0.1cm] {$\dots$}  
		-- (11, 1.6);
		\draw[blue!100] (8.6, 0.8) node[blue, right=0.7cm, above=0.1cm] {$2,k-1$}  
		node[blue, right=2.7cm, above=0.1cm] {$\dots$} -- (10, 0.8);
		\draw[blue!100] (8.6, -0.2) node[blue, right=1.5cm, above=0.3cm] {$\dots$}   
		-- (11.8, -0.2);
		\draw[blue!100] (8.6, -1) node[blue, right=0.7cm, above=0.1cm] {$2,4$} 
		node[blue, right=1.9cm, above=0.1cm] {$3,4$}  -- (11, -1);
		\draw[blue!100] (8.6, -1.8) node[blue, right=0.7cm, above=0.1cm] {$2,3$}   
		-- (10, -1.8);
		\draw[red, very thick] (-0.7,-1.8) rectangle (0.7, 3.2);
		\draw[red, very thick] (8.6,2.4) rectangle (14.6, 3.2);
		\draw[teal, very thick] (-0.7,3.2) rectangle (0.7, 4);
		\draw[teal, very thick] (14.6,3.2) rectangle (15.4, 4);
		
		\draw[->] (6.7,1.8) node[right=0.7cm, above] {$\pi$} -- (8,1.8);
		\node at (-1,1.1) {$k$};
		\node at (3.5,4.4) {$k$};
		\node at (8.2,1.1) {$k$};
		\node at (12.5,4.4) {$k$};
		
	\end{tikzpicture}
	\caption{Enumeration \ref{int:eq10} and the action of $\pi$}
	\label{fig:pi}
\end{figure}

%% file: pic_7.tex
\begin{figure}[t!]
	\centering
	\begin{tikzpicture}	
		\fill[color=blue!10] (0,3.2) rectangle (4,4);
		\draw[blue!100] (0, 0)  -- (0, 4);
		\draw[blue!100] (0, 4)  -- (4, 4);
		\draw[blue!100] (0, 3.2) node[blue, right=0.4cm, above=0.1cm] {$1,1$}  node[blue, right=1.2cm, above=0.1cm] {$2,2$} node[blue, right=2cm, above=0.1cm] {$3,3$} node[blue, right=2.8cm, above=0.1cm] {$4,4$} node[blue, right=3.6cm, above=0.1cm] {$5,5$} -- (4, 3.2);
		\draw[blue!100] (0.8, 0)  -- (0.8, 4);
		\draw[blue!100] (1.6, 0.8)  -- (1.6, 4);
		\draw[blue!100] (2.4, 1.6)  -- (2.4, 4);
		\draw[blue!100] (3.2, 2.4)  -- (3.2, 4);
		\draw[blue!100] (4, 3.2)  -- (4, 4);
		\draw[blue!100] (0, 2.4) node[blue, right=0.4cm, above=0.1cm] {$1,5$}  node[blue, right=1.2cm, above=0.1cm] {$2,5$} node[blue, right=2cm, above=0.1cm] {$3,5$} node[blue, right=2.8cm, above=0.1cm] {$4,5$}  -- (3.2, 2.4);
		\draw[blue!100] (0, 1.6) node[blue, right=0.4cm, above=0.1cm] {$1,4$}  node[blue, right=1.2cm, above=0.1cm] {$2,4$} node[blue, right=2cm, above=0.1cm] {$3,4$}  -- (2.4, 1.6);
		\draw[blue!100] (0, 0.8) node[blue, right=0.4cm, above=0.1cm] {$1,3$}  node[blue, right=1.2cm, above=0.1cm] {$2,3$} -- (1.6, 0.8);
		\draw[blue!100] (0, 0) node[blue, right=0.4cm, above=0.1cm] {$1,2$}   -- (0.8, 0);
		
		\fill[color=blue!10] (7,3.2) rectangle (11,4);
		\draw[blue!100] (7, 0)  -- (7, 4);
		\draw[blue!100] (7, 4)  -- (11, 4);
		\draw[blue!100] (7, 3.2) node[blue, right=0.4cm, above=0.1cm] {$2,2$}  
		node[blue, right=1.2cm, above=0.1cm] {$3,3$} 
		node[blue, right=2cm, above=0.1cm] {$4,4$} 
		node[blue, right=2.8cm, above=0.1cm] {$5,5$} 
		node[blue, right=3.6cm, above=0.1cm] {$1,1$} -- (11, 3.2);
		\draw[blue!100] (7.8, 0)  -- (7.8, 4);
		\draw[blue!100] (8.6, 0.8)  -- (8.6, 4);
		\draw[blue!100] (9.4, 1.6)  -- (9.4, 4);
		\draw[blue!100] (10.2, 2.4)  -- (10.2, 4);
		\draw[blue!100] (11, 3.2)  -- (11, 4);
		\draw[blue!100] (7, 2.4) node[blue, right=0.4cm, above=0.1cm] {$1,2$}  
		node[blue, right=1.2cm, above=0.1cm] {$1,3$} 
		node[blue, right=2cm, above=0.1cm] {$1,4$} 
		node[blue, right=2.8cm, above=0.1cm] {$1,5$}  -- (10.2, 2.4);
		\draw[blue!100] (7, 1.6) node[blue, right=0.4cm, above=0.1cm] {$2,5$}  
		node[blue, right=1.2cm, above=0.1cm] {$3,5$} 
		node[blue, right=2cm, above=0.1cm] {$4,5$}  -- (9.4, 1.6);
		\draw[blue!100] (7, 0.8) node[blue, right=0.4cm, above=0.1cm] {$2,4$}  
		node[blue, right=1.2cm, above=0.1cm] {$3,4$} -- (8.6, 0.8);
		\draw[blue!100] (7, 0) node[blue, right=0.4cm, above=0.1cm] {$2,3$}   -- (7.8, 0);
		\draw[->] (5,2) node[right=0.5cm, above] {$\pi$} -- (6,2);
		\draw[red, very thick] (0,0) rectangle (0.8, 3.2);
		\draw[red, very thick] (7,2.4) rectangle (10.2, 3.2);
		\draw[teal, very thick] (0,3.2) rectangle (0.8, 4);
		\draw[teal, very thick] (10.2,3.2) rectangle (11, 4);
	\end{tikzpicture}
\caption{The action of $\pi$ for $k = 5$ and $n = 10$}
\label{fig:perm}
\end{figure}
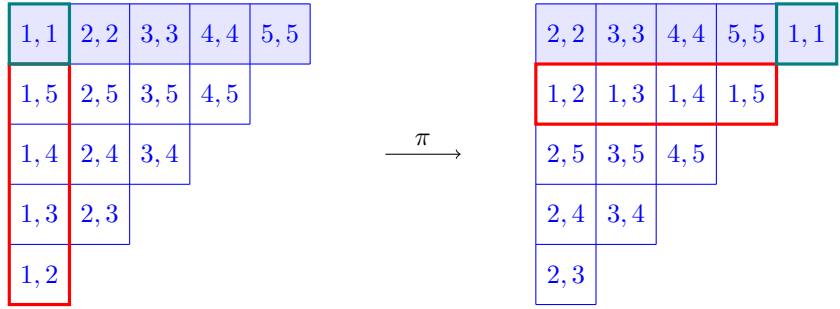

%% file: pic_6.tex
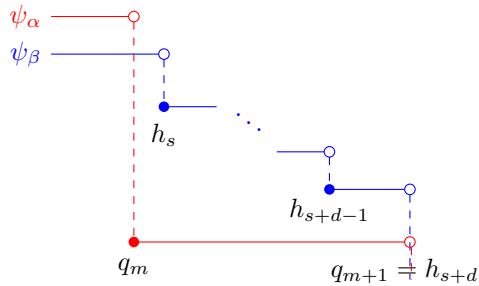
\begin{figure}[h]
	\centering
	\begin{tikzpicture}
		\draw[blue] (0,0.7) -- (-1, 0.7) node[black, below] {$h_{s+d-1}$};
		\draw[blue] (0.07,0.7) circle (0.07);
		\draw[blue, dashed] (0.07, 0.63) -- (0.07, -0.5);
		\fill[blue] (-1,0.7) circle (0.07);
		\draw[blue] (-1, 1.2) -- (-1.7, 1.2);
		\fill[white] (-1,1.2) circle (0.07);
		\draw[blue] (-1,1.2) circle (0.07);
		\draw[blue, dashed] (-1, 1.2-0.07) -- (-1, 0.7);
		\draw[blue] (-2.5, 1.8) node[ below=0.1cm, right=0.1cm] {$\ddots$} -- (-3.2, 1.8) node[black, below=0.1cm] {$h_s$};
		\fill[blue] (-3.2,1.8) circle (0.07);
		\draw[blue] (-3.2, 2.5) -- (-4.7, 2.5) node[left] {$\psi_{\beta}$};
		\fill[white] (-3.2,2.5) circle (0.07);
		\draw[blue] (-3.2,2.5) circle (0.07);
		\draw[blue, dashed] (-3.2, 2.5-0.07) -- (-3.2, 1.8);
		\draw[red] (0, 0) node[black, below=0.1cm] {$q_{m+1} = h_{s+d}$}-- (-3.6, 0) node[black, below=0.1cm] {$q_m$};
		\draw[red] (0.07,0) circle (0.07);
		\draw[red, dashed] (0.09, -0.07) -- (0.09, -0.5);
		\fill[red] (-3.6,0) circle (0.07);
		\draw[red] (-3.6, 3) -- (-4.7, 3) node[left] {$\psi_{\alpha}$};
		\fill[white] (-3.6,3) circle (0.07);
		\draw[red] (-3.6,3) circle (0.07);
		\draw[red, dashed] (-3.6, 3-0.07) -- (-3.6, 0);
	\end{tikzpicture}
	\caption{To Lemma \ref{lm:MM}}
	\label{fig1}
\end{figure}

%% file: pic_1.tex
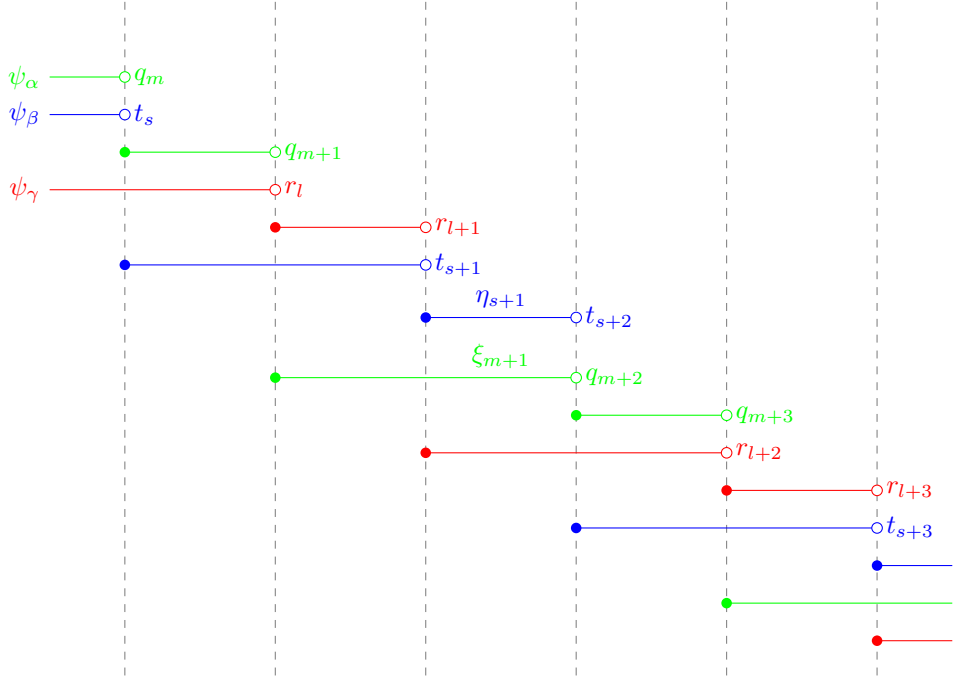
\begin{figure}
	\centering
	\begin{tikzpicture}
		\draw[gray!100, dashed] (0, 8)   -- (0, -1);
		\draw[gray!100, dashed] (2, 8)   -- (2, -1);
		\draw[gray!100, dashed] (4, 8)   -- (4, -1);
		\draw[gray!100, dashed] (6, 8)  -- (6, -1);
		\draw[gray!100, dashed] (8, 8)  -- (8, -1);
		\draw[gray!100, dashed] (10, 8)  -- (10, -1);
		\draw[green!100] (-1, 7) node[left] {$\psi_\alpha$} -- (0, 7) node[right] {$q_m$};
			\fill[white] (0,7) circle (0.07);
			\draw[green] (0,7) circle (0.07);
		\draw[blue!100] (-1, 6.5)node[left] {$\psi_\beta$} -- (0, 6.5) node[right] {$t_s$};
			\fill[white] (0,6.5) circle (0.07);
			\draw[blue] (0,6.5) circle (0.07);
		\draw[green!100] (0, 6) -- (2, 6)  node[right]  {$q_{m+1}$};
			\fill[green] (0,6) circle (0.07);
			\fill[white] (2,6) circle (0.07);
			\draw[green] (2,6) circle (0.07);
		\draw[red!100] (-1, 5.5) node[left] {$\psi_\gamma$} -- (2,5.5) node[right] {$r_l$};
			\fill[white] (2,5.5) circle (0.07);
			\draw[red] (2,5.5) circle (0.07);
		\draw[red!100] (2, 5) -- (4,5) node[right] {$r_{l+1}$};
			\fill[red] (2,5) circle (0.07);
			\fill[white] (4,5) circle (0.07);
			\draw[red] (4,5) circle (0.07);
		\draw[blue!100] (0, 4.5) -- (4, 4.5) node[right] {$t_{s+1}$};
			\fill[blue] (0,4.5) circle (0.07);
			\fill[white] (4,4.5) circle (0.07);
			\draw[blue] (4,4.5) circle (0.07);
		\draw[green!100] (2, 3) -- (6,3) node[left=1cm, above] {$\xi_{m+1}$} node[right] {$q_{m+2}$};
			\fill[green] (2,3) circle (0.07);
			\fill[white] (6,3) circle (0.07);
			\draw[green] (6,3) circle (0.07);
		\draw[blue!100] (4, 3.8) -- (6, 3.8) node[midway, above] {$\eta_{s+1}$} node[right] {$t_{s+2}$};
			\fill[blue] (4,3.8) circle (0.07);
			\fill[white] (6,3.8) circle (0.07);
			\draw[blue] (6,3.8) circle (0.07);
		\draw[green!100] (6, 2.5) -- (8, 2.5) node[right] {$q_{m+3}$};
			\fill[green] (6,2.5) circle (0.07);
			\fill[white] (8,2.5) circle (0.07);
			\draw[green] (8,2.5) circle (0.07);
		\draw[red!100] (4, 2) -- (8, 2) node[right] {$r_{l+2}$};
			\fill[red] (4,2) circle (0.07);
			\fill[white] (8,2) circle (0.07);
			\draw[red] (8,2) circle (0.07);
		\draw[red!100] (8, 1.5) -- (10, 1.5) node[right] {$r_{l+3}$};
			\fill[red] (8,1.5) circle (0.07);
			\fill[white] (10,1.5) circle (0.07);
			\draw[red] (10,1.5) circle (0.07);
		\draw[blue!100] (6, 1) -- (10, 1) node[right] {$t_{s+3}$};
			\fill[blue] (6,1) circle (0.07);
			\fill[white] (10,1) circle (0.07);
			\draw[blue] (10,1) circle (0.07);
		\draw[green!100] (8, 0) -- (11, 0);
			\fill[green] (8,0) circle (0.07);
		\draw[blue!100] (10, 0.5) -- (11, 0.5);
			\fill[blue] (10,0.5) circle (0.07);
		\draw[red!100] (10, -0.5) -- (11, -0.5);
			\fill[red] (10, -0.5) circle (0.07);
	\end{tikzpicture}
	\caption{to Lemma \ref{lm:3_func}}
	\label{fig_lm_3func}
\end{figure}

%% file: pic_11.tex
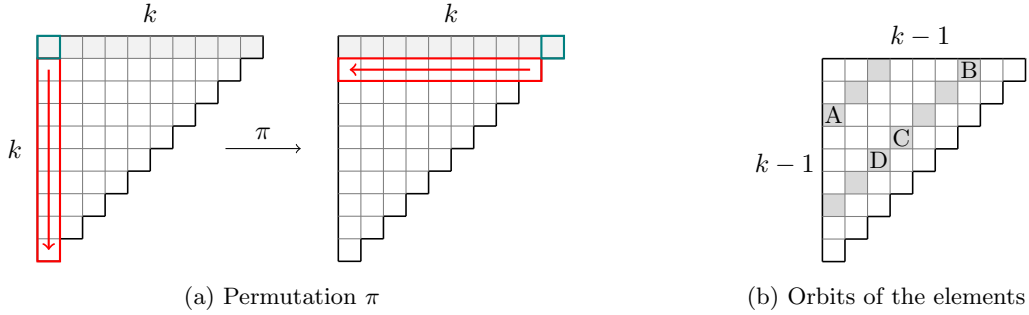
\begin{figure}[t!]
	\centering
	\begin{subfigure}{.45\linewidth}
		\centering
		\begin{tikzpicture}	
			\fill[color=gray!10] (0,2.7) rectangle (3,3);
			\fill[color=gray!10] (4,2.7) rectangle (7,3);
			\draw[semithick] (0,0) -- (0,3);
			\draw[semithick] (0,3) -- (3,3);
			\draw[semithick] (0,0) -- (0.3,0);
			\draw[semithick] (0.3,0) -- (0.3,0.3);
			\draw[semithick] (0.3,0.3) -- (0.6,0.3);
			\draw[semithick] (0.6,0.3) -- (0.6,0.6);
			\draw[semithick] (0.6,0.6) -- (0.9,0.6);
			\draw[semithick] (0.9,0.6) -- (0.9,0.9);
			\draw[semithick] (0.9,0.9) -- (1.2,0.9);
			\draw[semithick] (1.2,0.9) -- (1.2,1.2);
			\draw[semithick] (1.2,1.2) -- (1.5,1.2);
			\draw[semithick] (1.5,1.2) -- (1.5,1.5);
			\draw[semithick] (1.5,1.5) -- (1.8,1.5);
			\draw[semithick] (1.8,1.5) -- (1.8,1.8);
			\draw[semithick] (1.8,1.8) -- (2.1,1.8);
			\draw[semithick] (2.1,1.8) -- (2.1,2.1);
			\draw[semithick] (2.1,2.1) -- (2.4,2.1);
			\draw[semithick] (2.4,2.1) -- (2.4,2.4);
			\draw[semithick] (2.4,2.4) -- (2.7,2.4);
			\draw[semithick] (2.7,2.4) -- (2.7,2.7);
			\draw[semithick] (2.7,2.7) -- (3.0,2.7);
			\draw[semithick] (3.0,2.7) -- (3.0,3.0);
			
			\draw[gray] (0.3,0.3) -- (0.3,3);
			\draw[gray] (0.6,0.6) -- (0.6,3);
			\draw[gray] (0.9,0.9) -- (0.9,3);
			\draw[gray] (1.2,1.2) -- (1.2,3);
			\draw[gray] (1.5,1.5) -- (1.5,3);
			\draw[gray] (1.8,1.8) -- (1.8,3);
			\draw[gray] (2.1,2.1) -- (2.1,3);
			\draw[gray] (2.4,2.4) -- (2.4,3);
			\draw[gray] (2.7,2.7) -- (2.7,3);
			\draw[gray] (3.0,3.0) -- (3.0,3);
			
			\draw[gray] (0,0.3) -- (0.3,0.3);
			\draw[gray] (0,0.6) -- (0.6,0.6);
			\draw[gray] (0,0.9) -- (0.9,0.9);
			\draw[gray] (0,1.2) -- (1.2,1.2);
			\draw[gray] (0,1.5) -- (1.5,1.5);
			\draw[gray] (0,1.8) -- (1.8,1.8);
			\draw[gray] (0,2.1) -- (2.1,2.1);
			\draw[gray] (0,2.4) -- (2.4,2.4);
			\draw[gray] (0,2.7) -- (2.7,2.7);
			\draw[gray] (0,3.0) -- (3.0,3.0);
			
			\draw[semithick] (4,0) -- (4,3);
			\draw[semithick] (4,3) -- (7,3);
			\draw[semithick] (4,0) -- (4.3,0);
			\draw[semithick] (4.3,0) -- (4.3,0.3);
			\draw[semithick] (4.3,0.3) -- (4.6,0.3);
			\draw[semithick] (4.6,0.3) -- (4.6,0.6);
			\draw[semithick] (4.6,0.6) -- (4.9,0.6);
			\draw[semithick] (4.9,0.6) -- (4.9,0.9);
			\draw[semithick] (4.9,0.9) -- (5.2,0.9);
			\draw[semithick] (5.2,0.9) -- (5.2,1.2);
			\draw[semithick] (5.2,1.2) -- (5.5,1.2);
			\draw[semithick] (5.5,1.2) -- (5.5,1.5);
			\draw[semithick] (5.5,1.5) -- (5.8,1.5);
			\draw[semithick] (5.8,1.5) -- (5.8,1.8);
			\draw[semithick] (5.8,1.8) -- (6.1,1.8);
			\draw[semithick] (6.1,1.8) -- (6.1,2.1);
			\draw[semithick] (6.1,2.1) -- (6.4,2.1);
			\draw[semithick] (6.4,2.1) -- (6.4,2.4);
			\draw[semithick] (6.4,2.4) -- (6.7,2.4);
			\draw[semithick] (6.7,2.4) -- (6.7,2.7);
			\draw[semithick] (6.7,2.7) -- (7.0,2.7);
			\draw[semithick] (7.0,2.7) -- (7.0,3.0);
			
			\draw[gray] (4.3,0.3) -- (4.3,3);
			\draw[gray] (4.6,0.6) -- (4.6,3);
			\draw[gray] (4.9,0.9) -- (4.9,3);
			\draw[gray] (5.2,1.2) -- (5.2,3);
			\draw[gray] (5.5,1.5) -- (5.5,3);
			\draw[gray] (5.8,1.8) -- (5.8,3);
			\draw[gray] (6.1,2.1) -- (6.1,3);
			\draw[gray] (6.4,2.4) -- (6.4,3);
			\draw[gray] (6.7,2.7) -- (6.7,3);
			\draw[gray] (7.0,3.0) -- (7.0,3);
			
			\draw[gray] (4,0.3) -- (4.3,0.3);
			\draw[gray] (4,0.6) -- (4.6,0.6);
			\draw[gray] (4,0.9) -- (4.9,0.9);
			\draw[gray] (4,1.2) -- (5.2,1.2);
			\draw[gray] (4,1.5) -- (5.5,1.5);
			\draw[gray] (4,1.8) -- (5.8,1.8);
			\draw[gray] (4,2.1) -- (6.1,2.1);
			\draw[gray] (4,2.4) -- (6.4,2.4);
			\draw[gray] (4,2.7) -- (6.7,2.7);
			\draw[gray] (4,3.0) -- (7.0,3.0);
			
			\draw[->] (2.5,1.5) node[right=0.5cm, above] {$\pi$} -- (3.5,1.5);
			
			\draw[red, thick] (0,0) rectangle (0.3, 2.7);
			\draw[red, thick] (4,2.4) rectangle (6.7, 2.7);
			\draw[teal, thick] (0,2.7) rectangle (0.3, 3);
			\draw[teal, thick] (6.7,2.7) rectangle (7, 3);
			\draw[->, red, thick] (0.15,2.55) -- (0.15,0.15);
			\draw[->, red, thick] (6.55,2.55) -- (4.15,2.55);
			\node at (-0.3,1.5) {$k$};
			\node at (1.5,3.3) {$k$};
			\node at (5.5,3.3) {$k$};
		\end{tikzpicture}
		\caption{Permutation $\pi$}
		\label{fig:pi:sub1}
	\end{subfigure}
	\begin{subfigure}{.45\linewidth}
		\centering
		\begin{tikzpicture}	
			\fill[color=gray!30] (0,1.8) rectangle (0.3,2.1) node[black, pos=.5] {\small A};	
			\fill[color=gray!30] (0.3,2.1) rectangle (0.6,2.4);
			\fill[color=gray!30] (0.6,2.4) rectangle (0.9,2.7);
			
			\fill[color=gray!30] (0,0.6) rectangle (0.3,0.9);	
			\fill[color=gray!30] (0.3,0.9) rectangle (0.6,1.2);
			\fill[color=gray!30] (0.6,1.2) rectangle (0.9,1.5) node[black, pos=.5] {\small D};
			\fill[color=gray!30] (0.9,1.5) rectangle (1.2,1.8) node[black, pos=.5] {\small C};
			\fill[color=gray!30] (1.2,1.8) rectangle (1.5,2.1);
			\fill[color=gray!30] (1.5,2.1) rectangle (1.8,2.4);
			\fill[color=gray!30] (1.8,2.4) rectangle (2.1,2.7) node[black, pos=.5] {\small B};

			\draw[semithick] (0,0) -- (0,2.7);
			\draw[semithick]  (0,2.7) -- (2.7,2.7);
			\draw[semithick] (0,0) -- (0.3,0);
			\draw[semithick] (0.3,0) -- (0.3,0.3);
			\draw[semithick] (0.3,0.3) -- (0.6,0.3);
			\draw[semithick] (0.6,0.3) -- (0.6,0.6);
			\draw[semithick] (0.6,0.6) -- (0.9,0.6);
			\draw[semithick] (0.9,0.6) -- (0.9,0.9);
			\draw[semithick] (0.9,0.9) -- (1.2,0.9);
			\draw[semithick] (1.2,0.9) -- (1.2,1.2);
			\draw[semithick] (1.2,1.2) -- (1.5,1.2);
			\draw[semithick] (1.5,1.2) -- (1.5,1.5);
			\draw[semithick] (1.5,1.5) -- (1.8,1.5);
			\draw[semithick] (1.8,1.5) -- (1.8,1.8);
			\draw[semithick] (1.8,1.8) -- (2.1,1.8);
			\draw[semithick] (2.1,1.8) -- (2.1,2.1);
			\draw[semithick] (2.1,2.1) -- (2.4,2.1);
			\draw[semithick] (2.4,2.1) -- (2.4,2.4);
			\draw[semithick] (2.4,2.4) -- (2.7,2.4);
			\draw[semithick] (2.7,2.4) -- (2.7,2.7);
			
			\draw[gray] (0.3,0.3) -- (0.3,2.7);
			\draw[gray] (0.6,0.6) -- (0.6,2.7);
			\draw[gray] (0.9,0.9) -- (0.9,2.7);
			\draw[gray] (1.2,1.2) -- (1.2,2.7);
			\draw[gray] (1.5,1.5) -- (1.5,2.7);
			\draw[gray] (1.8,1.8) -- (1.8,2.7);
			\draw[gray] (2.1,2.1) -- (2.1,2.7);
			\draw[gray] (2.4,2.4) -- (2.4,2.7);
			\draw[gray] (2.7,2.7) -- (2.7,2.7);

			\draw[gray] (0,0.3) -- (0.3,0.3);
			\draw[gray] (0,0.6) -- (0.6,0.6);
			\draw[gray] (0,0.9) -- (0.9,0.9);
			\draw[gray] (0,1.2) -- (1.2,1.2);
			\draw[gray] (0,1.5) -- (1.5,1.5);
			\draw[gray] (0,1.8) -- (1.8,1.8);
			\draw[gray] (0,2.1) -- (2.1,2.1);
			\draw[gray] (0,2.4) -- (2.4,2.4);
			
			\node at (-0.5,1.3) {$k-1$};
			\node at (1.3,3) {$k-1$};
			
		\end{tikzpicture}
		\caption{Orbits of the elements}
		\label{fig:pi:sub2}
	\end{subfigure}
	\caption{to Lemma \ref{lm:order}}
\end{figure}

%% file: pic_9.tex
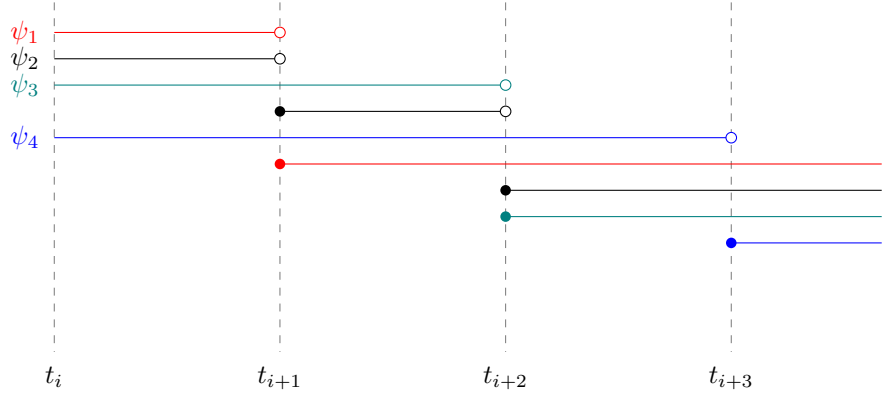
\begin{figure}[!t]
	\center
	\begin{tikzpicture}
		\draw[gray, dashed] (0,0.4)  -- (0,-4.25) node[black, below=2pt]{$t_i$};
		\draw[gray, dashed] (3,0.4)  -- (3,-4.25) node[black, below=2pt]{$t_{i+1}$};
		\draw[gray, dashed] (6,0.4)  -- (6,-4.25) node[black, below=2pt]{$t_{i+2}$};
		\draw[gray, dashed] (9,0.4)  -- (9,-4.25) node[black, below=2pt]{$t_{i+3}$};
		
		\draw[red] (0,0) node[left=2pt]{$\psi_1$} -- (3,0);
		\fill[white] (3,0) circle (0.07);
		\draw[red] (3,0) circle (0.07);
		\draw[black] (0,-0.35) node[left=2pt]{$\psi_2$} -- (3,-0.35);
		\fill[white] (3,-0.35) circle (0.07);
		\draw[black] (3,-0.35) circle (0.07);
		\draw[teal] (0,-0.7) node[left=2pt]{$\psi_3$} -- (6,-0.7);
		\fill[white] (6,-0.7) circle (0.07);
		\draw[teal] (6,-0.7) circle (0.07);
		
		\draw[black] (3,-1.05)  -- (6,-1.05);
		\fill[black] (3,-1.05) circle (0.07);
		\fill[white] (6,-1.05) circle (0.07);
		\draw[black] (6,-1.05) circle (0.07);
		
		\draw[blue] (0,-1.4) node[left=2pt]{$\psi_4$} -- (9,-1.4);
		\fill[white] (9,-1.4) circle (0.07);
		\draw[blue] (9,-1.4) circle (0.07);
		
		\draw[red] (3,-1.75)  -- (11,-1.75);
		\fill[red] (3,-1.75) circle (0.07);
		\draw[black] (6,-2.1)  -- (11,-2.1);
		\fill[black] (6,-2.1) circle (0.07);
		
		\draw[teal] (6,-2.45)  -- (11,-2.45);
		\fill[teal] (6,-2.45) circle (0.07);
		
		\draw[blue] (9,-2.8)  -- (11,-2.8);
		\fill[blue] (9,-2.8) circle (0.07);
	\end{tikzpicture}
	\caption{The example for projection operator $\pr_{\bbeta}$. Picture continues periodically.}
	\label{fig:pr}
\end{figure}

%% file: pic_3.tex
\begin{figure}[b!]
	\centering
	\begin{tikzpicture}	
		\draw[gray!100, dashed] (0, 7.3)  -- (0,0.5) node[black, left] {$t_i$};
		\draw[gray!100, dashed] (2.5, 7.3) node[black, left=1cm] {$\blv_1$} node[black, right=1cm] {$\blv_2$}  -- (2.5,0.5) node[black, left] {$t_{i+1}$};
		\draw[gray!100, dashed] (5, 7.3) node[black, right=1cm] {$\blv_1$}  -- (5,0.5) node[black, left] {$t_{i+2}$};
		\draw[gray!100, dashed] (7.5, 7.3) node[black, right=1cm] {$\blv_3$}  -- (7.5,0.5) node[black, left] {$t_{i+3}$};
		\draw[gray!100, dashed] (10, 7.3) node[black, right=1cm] {$\blv_{1/2}$}  -- (10,0.5) node[black, left] {$t_{i+4}$};
		%\draw[red!100] (0, 7.65) node[left] {$\psi_{\a_4}$} -- (2.5, 7.65);
		%\fill[white] (2.5,7.65) circle (0.05);
		%\draw[red] (2.5,7.65) circle (0.05);
		\draw[red!100] (0, 6.6) node[left] {$\psi_{1}$} -- (7.5, 6.6);
		\fill[white] (7.5,6.6) circle (0.05);
		\draw[red] (7.5,6.6) circle (0.05);
		%\draw[red!100] (2.5, 6.95)  -- (7.5, 6.95);
		%	\fill[red] (2.5,6.95) circle (0.05);
		%\fill[white] (7.5,6.95) circle (0.05);
		%\draw[red] (7.5,6.95) circle (0.05);
		%\draw[orange!100] (5, 6.6)  -- (7.5, 6.6);
		%	\fill[orange] (5,6.6) circle (0.05);
		%	\fill[white] (7.5,6.6) circle (0.05);
		%	\draw[orange] (7.5,6.6) circle (0.05);
		\draw[cyan!100] (0, 6.25) node[left] {$\psi_{2}$}  -- (7.5, 6.25);
		\fill[white] (7.5,6.25) circle (0.05);
		\draw[cyan] (7.5,6.25) circle (0.05);
		\draw[green!100] (0,5.9) node[left] {$\psi_{4}$} -- (2.5, 5.9);
		\fill[white] (2.5,5.9) circle (0.05);
		\draw[green] (2.5,5.9) circle (0.05);
		\draw[blue!100] (0,5.55) node[left] {$\psi_{5}$} -- (2.5, 5.55);
		\fill[white] (2.5,5.55) circle (0.05);
		\draw[blue] (2.5,5.55) circle (0.05);
		\draw[black!100] (0,5.2)  node[left] {$\psi_{3}$} -- (5, 5.2);
		\fill[white] (5,5.2) circle (0.05);
		\draw[black] (5,5.2) circle (0.05);
		\draw[green!100] (2.5,4.5)  -- (10, 4.5);
		\fill[green] (2.5,4.5) circle (0.05);
		\fill[white] (10,4.5) circle (0.05);
		\draw[green] (10,4.5) circle (0.05);
		\draw[blue!100] (2.5,4.85)  -- (5, 4.85);
		\fill[blue] (2.5,4.85) circle (0.05);
		\fill[white] (5,4.85) circle (0.05);
		\draw[blue] (5,4.85) circle (0.05);
		\draw[blue!100] (5,4.15)  -- (10, 4.15);
		\fill[blue] (5,4.15) circle (0.05);
		\fill[white] (10,4.15) circle (0.05);
		\draw[blue] (10,4.15) circle (0.05);
		\draw[black!100] (5,3.8)  -- (10, 3.8);
		\fill[black] (5,3.8) circle (0.05);
		\fill[white] (10,3.8) circle (0.05);
		\draw[black] (10,3.8) circle (0.05);
		\draw[cyan!100] (7.5, 3.45)  -- (10, 3.45);
		\fill[cyan] (7.5,3.45) circle (0.05);
		\fill[white] (10,3.45) circle (0.05);
		\draw[cyan] (10,3.45) circle (0.05);
		%\draw[red] (7.5, 3.1) -- (12.5, 3.1);
		%\fill[red] (7.5,3.1) circle (0.05);
		\draw[red!100] (7.5, 2.75)  -- (12.5, 2.75);
		\fill[red] (7.5,2.75) circle (0.05);
		\draw[cyan!100] (10, 2.4)  -- (12.5, 2.4);
		\fill[cyan] (10,2.4) circle (0.05);
		\draw[green!100] (10, 2.05)  -- (12.5, 2.05);
		\fill[green] (10,2.05) circle (0.05);
		\draw[blue!100] (10, 1.7)  -- (12.5, 1.7);
		\fill[blue] (10,1.7) circle (0.05);
		\draw[black!100] (10, 1.35)  -- (12.5, 1.35);
		\fill[black] (10,1.35) circle (0.05);
	\end{tikzpicture}
	\caption{Base of induction for the proof of \ref{it2} from Theorem \ref{th:main}}
	\label{fig:period}
\end{figure}
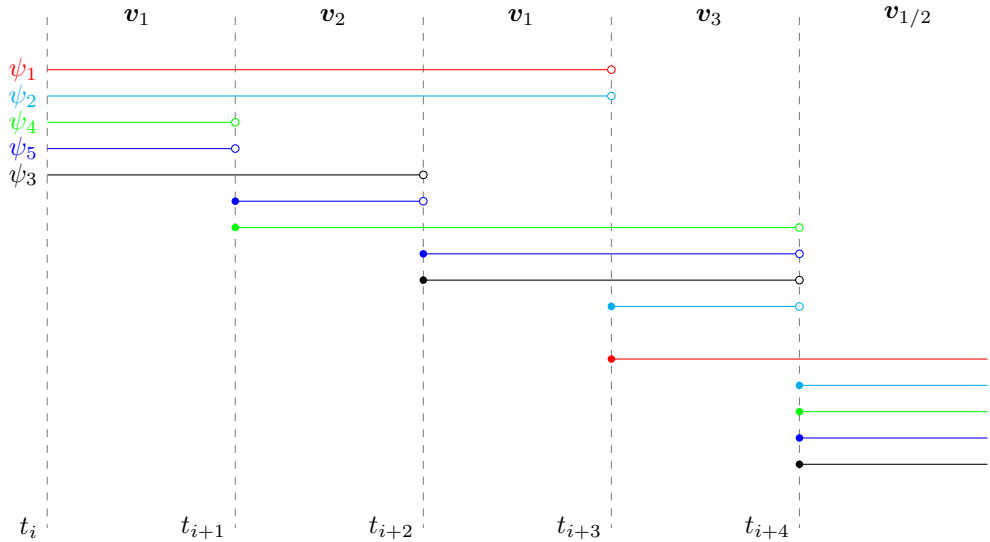

%% file: pic_8.tex
\begin{figure}[t]
	\centering
	\begin{tikzpicture}	
		\draw[gray, dashed] (0, 0) node[teal, above=0.1cm, right=0.2cm] {$\blw_1$}  node[blue, above=0.7cm, right=0.2cm] {$\blu_1$} node[black, above=1.2cm, right=0.2cm] {$\blv_1$} node[black, above=1.8cm, right=0.2cm] {$t_0$} -- (0, 1.5) ;
		\draw[gray, dashed] (1.3, 0) node[teal, above=0.1cm, right=0.2cm] {$\blw_2$}  node[blue, above=0.7cm, right=0.2cm] {$\blu_1$} node[black, above=1.2cm, right=0.2cm] {$\blv_2$} node[black, above=1.8cm, right=0.2cm] {$t_1(\a)$} -- (1.3, 1.5) ;
		\draw[gray, dashed] (2.6, 0) node[teal, above=0.1cm, right=0.2cm] {$\blw_2$}  node[blue, above=0.7cm, right=0.2cm] {$\blu_1$} node[black, above=1.2cm, right=0.2cm] {$\blv_1$} node[black, above=1.8cm, right=0.2cm] {$t_2(\a)$} -- (2.6, 1.5) ;
		\draw[gray, dashed] (3.9, 0) node[black, above=0.7cm, right=0.7cm] {$\dots$} -- (3.9, 1.5) ;
		\draw[gray, dashed] (6.2, 0) node[teal, above=0.1cm, right=0.5cm] {$\blw_{2}$}  node[blue, above=0.7cm, right=0.5cm] {$\blu_1$} node[black, above=1.2cm, right=0.5cm] {$\blv_{1/2}$}  node[black, above=1.8cm, right=0.2cm] {$t_{s+1}(\a)$}-- (6.2, 1.5);
		\draw[gray, dashed] (8, 0) node[teal, above=0.1cm, right=1.2cm] {$\blw_{2}$}  node[blue, above=0.7cm, right=1.2cm] {$\blu_2$} node[black, above=1.2cm, right=1.2cm] {$\blv_{3}$} node[black, above=1.8cm, right=0.2cm] {$t_{s+2}(\a) = t_1(\bbeta_1)$} -- (8, 1.5) ;
		\draw[gray, dashed] (11, 0) node[teal, above=0.1cm, right=1cm] {$\blw_{3}$}  node[blue, above=0.7cm, right=1cm] {$\blu_3$} node[black, above=1.2cm, right=1cm] {$\blv_{4}$} node[black, above=1.8cm, right=0.2cm] {$t_{s+3}(\a) = t_2(\bbeta_1)$} node[black, above=0.7cm, right=3cm] {$\dots$}-- (11, 1.5) ;
	\end{tikzpicture}
	\caption{Step of induction for \ref{it2}}
	\label{fig:per}
\end{figure}

%% file: pic_12.tex
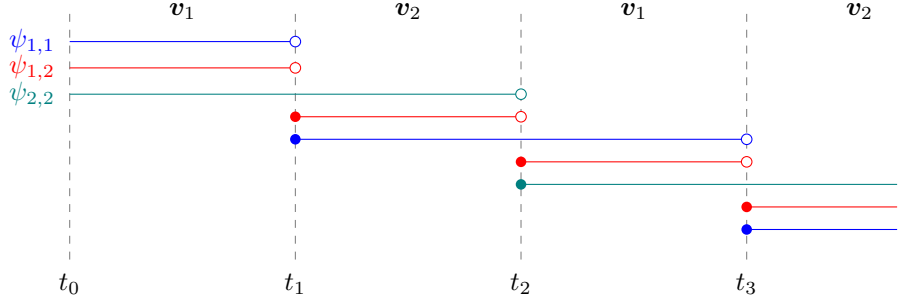
\begin{figure}[!t]
	\center
	\begin{tikzpicture}
		\draw[gray, dashed] (0,0) node[black, below=2]{$t_0$} -- (0,3.3) node[black,right=1.2cm] {$\blv_1$};
		\draw[gray, dashed] (3,0) node[black, below=2]{$t_1$} -- (3,3.3) node[black,right=1.2cm] {$\blv_2$};
		\draw[gray, dashed] (6,0) node[black, below=2]{$t_2$} -- (6,3.3) node[black,right=1.2cm] {$\blv_1$};
		\draw[gray, dashed] (9,0) node[black, below=2]{$t_3$} -- (9,3.3) node[black,right=1.2cm] {$\blv_2$};
		\draw[blue] (0, 2.9) node[left=2]{$\psi_{1,1}$} -- (3, 2.9);
		\fill[white] (3,2.9) circle (0.07);
		\draw[blue] (3,2.9) circle (0.07);
		\draw[red] (0, 2.55) node[left=2]{$\psi_{1,2}$} -- (3, 2.55);
		\fill[white] (3,2.55) circle (0.07);
		\draw[red] (3,2.55) circle (0.07);
		\draw[teal] (0, 2.2) node[left=2]{$\psi_{2,2}$} -- (6, 2.2);
		\fill[white] (6,2.2) circle (0.07);
		\draw[teal] (6,2.2) circle (0.07);
		\draw[red] (3, 1.9) -- (6, 1.9);
		\fill[red] (3,1.9) circle (0.07);
		\fill[white] (6,1.9) circle (0.07);
		\draw[red] (6,1.9) circle (0.07);
		\draw[blue] (3, 1.6) -- (9, 1.6);
		\fill[blue] (3,1.6) circle (0.07);
		\fill[white] (9,1.6) circle (0.07);
		\draw[blue] (9,1.6) circle (0.07);
		\draw[red] (6, 1.3) -- (9, 1.3);
		\fill[red] (6,1.3) circle (0.07);
		\fill[white] (9,1.3) circle (0.07);
		\draw[red] (9,1.3) circle (0.07);
		\draw[teal] (6, 1) -- (11, 1);
		\fill[teal] (6,1) circle (0.07);
		\draw[red] (9, 0.7) -- (11, 0.7);
		\fill[red] (9,0.7) circle (0.07);
		\draw[blue] (9, 0.4) -- (11, 0.4);
		\fill[blue] (9,0.4) circle (0.07);
		
	\end{tikzpicture}
	\caption{Case $k=2$ and $n=3$}
	\label{fig:ex3}
\end{figure}

%% file: pic_2.tex
\begin{figure}[!t]
	\center
	\begin{tikzpicture}
		\draw[gray, dashed] (-2,1) node[black,right=1.2cm] {$\blv_1$} node[teal,below=0.5cm, right=1.2cm] {$\blu_1$}-- (-2,-6.9) node[black, below=2pt]{$t_0$};
		\draw[gray, dashed] (1,1) node[black, right=1.2cm] {$\blv_2$} node[teal,below=0.5cm, right=1.2cm] {$\blu_2$} -- (1,-6.9) node[black, below=2pt]{$t_1(\a) = t_1(\bbeta_1)$};
		\draw[gray, dashed] (4,1) node[black, right=1.2cm] {$\blv_3$} node[teal,below=0.5cm, right=1.2cm] {$\blu_1$} -- (4,-6.9) node[black, below=2pt]{$t_2(\a) = t_2(\bbeta_1)$};
		\draw[gray, dashed] (7,1) node[black, right=1.2cm] {$\blv_1$} node[teal,below=0.5cm, right=1.2cm] {$\blu_1$} -- (7,-6.9) node[black, below=2pt]{$t_3(\a)$};
		\draw[gray, dashed] (10,1) node[black, right=1.2cm] {$\blv_2$} node[teal,below=0.5cm, right=1.2cm] {$\blu_2$} -- (10,-6.9) node[black, below=2pt]{$t_4(\a)=t_3(\bbeta_1)$};
		
		\draw[red] (-2, 0) node[left=2pt]{$\psi_1 =\psi_{1,1}$} -- (1, 0);
		\fill[white] (1,0) circle (0.07);
		\draw[red] (1,0) circle (0.07);
		
		\draw[yellow!40!orange] (-2, -0.35) node[left=2pt]{$\psi_2 = \psi_{1,3}$} -- (1, -0.35);
		\fill[white] (1,-0.35) circle (0.07);
		\draw[yellow!40!orange] (1,-0.35) circle (0.07);
		
		\draw[gray] (-2, -0.7) node[left=2pt]{$\psi_3= \psi_{1,2}$} -- (1, -0.7);
		\fill[white] (1,-0.7) circle (0.07);
		\draw[gray] (1,-0.7) circle (0.07);
		
		\draw[green] (-2, -1.2) node[left=2pt]{$\psi_4 =\psi_{2,2}$} -- (4, -1.2);
		\fill[white] (4,-1.2) circle (0.07);
		\draw[green] (4,-1.2) circle (0.07);
		
		\draw[blue] (-2, -1.8) node[left=2pt]{$\psi_5 =\psi_{2,3}$} -- (4, -1.8);
		\fill[white] (4,-1.8) circle (0.07);
		\draw[blue] (4,-1.8) circle (0.07);
		
		\fill[gray] (1,-1.5) circle (0.07);
		\draw[gray] (1, -1.5) -- (4, -1.5);
		\fill[white] (4,-1.5) circle (0.07);
		\draw[gray] (4,-1.5) circle (0.07);
		
		\draw[black] (-2,-2.3) node[left=2pt]{$\psi_6 =\psi_{3,3}$} -- (7, -2.3);
		\fill[white] (7,-2.3) circle (0.07);
		\draw[black] (7,-2.3) circle (0.07);
		
		\fill[yellow!40!orange] (1,-2.9) circle (0.07);
		\draw[yellow!40!orange] (1, -2.9) -- (7, -2.9);
		\fill[white] (7,-2.9) circle (0.07);
		\draw[yellow!40!orange] (7,-2.9) circle (0.07);
		
		\fill[blue] (4,-2.6) circle (0.07);
		\draw[blue] (4, -2.6) -- (7, -2.6);
		\fill[white] (7,-2.6) circle (0.07);
		\draw[blue] (7,-2.6) circle (0.07);
		
		\draw[red] (1, -3.4) -- (10, -3.4);
		\fill[red] (1,-3.4) circle (0.07);
		\fill[white] (10,-3.4) circle (0.07);
		\draw[red] (10,-3.4) circle (0.07);
		
		\draw[yellow!40!orange] (7, -3.7) -- (10, -3.7);
		\fill[yellow!40!orange] (7,-3.7) circle (0.07);
		\fill[white] (10,-3.7) circle (0.07);
		\draw[yellow!40!orange] (10,-3.7) circle (0.07);
		
		\draw[gray] (4, -4) -- (10, -4);
		\fill[gray] (4,-4) circle (0.07);
		\fill[white] (10,-4) circle (0.07);
		\draw[gray] (10,-4) circle (0.07);
		
		\draw[green] (4, -4.5) -- (13, -4.5);
		\fill[green] (4,-4.5) circle (0.07);
		
		\fill[gray] (10,-4.8) circle (0.07);
		\draw[gray] (10, -4.8) -- (13, -4.8);
		
		\fill[blue] (7,-5.1) circle (0.07);
		\draw[blue] (7, -5.1) -- (13, -5.1);
		
		\fill (7,-5.6) circle (0.07);
		\draw (7,-5.6) -- (13, -5.6);
		
		\fill[yellow!40!orange] (10,-5.9) circle (0.07);
		\draw[yellow!40!orange] (10, -5.9) -- (13, -5.9);
		
		\fill[red] (10,-6.2) circle (0.07);
		\draw[red] (10, -6.2) -- (13, -6.2);
		
	\end{tikzpicture}
	\caption{Case $k=3, \ n=6$}
	\label{fig2}
\end{figure}
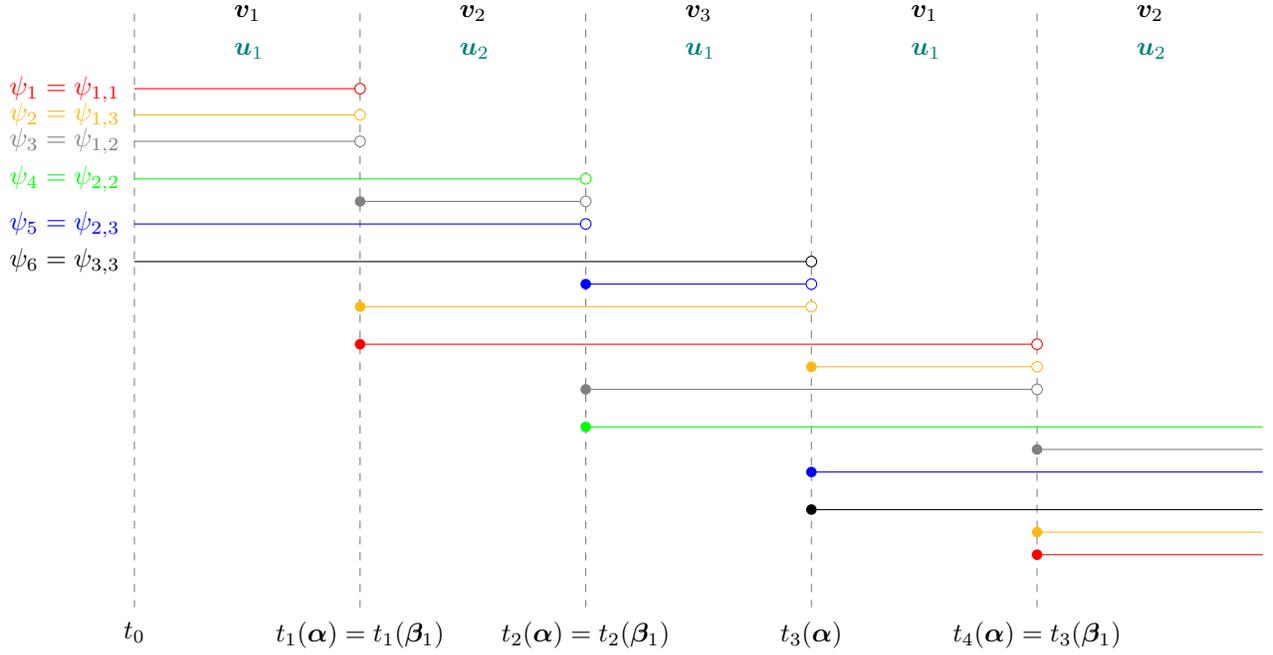

%% file: bibliography.tex
\bibliographystyle{abbrv}
%%%